\documentclass[12 pt,twosided,reqno]{amsart}
\usepackage[arrow,matrix]{xy}
\usepackage{graphicx}
\usepackage{epsfig}

\usepackage{amscd,amsmath,amsthm,amssymb}
\usepackage{pstcol,pst-plot,pst-3d}

\theoremstyle{plain}
\newtheorem{theorem}[subsection]{Theorem}
\newtheorem{lemma}[subsection]{Lemma}
\newtheorem{proposition}[subsection]{Proposition}
\newtheorem{cor}[subsection]{Corollary}
\theoremstyle{definition}
\newtheorem{remark}[subsection]{Remark}
\newtheorem{definition}[subsection]{Definition}
\newtheorem{example}[subsection]{Example}

\newcommand{\SQR}{\put(15,0){\line(-1,1){15}}
\put(15,0){\line(1,1){15}}
\put(15,30){\line(-1,-1){15}}
\put(15,30){\line(1,-1){15}}}

\def\Div{\hbox{Div}}

\begin{document}
\title [What Could be a  Simple Permutation?]
{What Could be a  Simple Permutation?}
\thanks{\emph{keywords and phrases:}  permutations, Coxeter generators,  Cayley graph}
\thanks{This research is partially supported by Higher Education Commission, Pakistan.\\
2010 AMS classification: Primary 05A05, 05A15, 05A16; Secondary 20B19, 20F36.}
\author{  REHANA ASHRAF$^{1}$,\,\,BARBU BERCEANU$^{1,2}$,\,\, AYESHA RIASAT$^{1}$}

\address{$^{1}$Abdus Salam School of Mathematical Sciences,
 GC University, Lahore-Pakistan.}
\email {rashraf@sms.edu.pk}
\email {ayesha.riyasat@gmail.com}
\address{$^{2}$
 Institute of Mathematics Simion Stoilow, Bucharest-Romania (permanent address).}
\email {Barbu.Berceanu@imar.ro}
  \maketitle
 \pagestyle{myheadings} \markboth{\centerline {\scriptsize
 REHANA ASHRAF,\,\,\,BARBU BERCEANU,\,\,\,AYESHA RIASAT   }} {\centerline {\scriptsize
  What Could be a  Simple Permutation?}}
\begin{abstract} Different ways to describe a permutation, as a  sequence of integers, or a  product of Coxeter generators, or  a tree, give different choices to define a simple permutation. We recollect  few of them, define new types of simple permutations, and analyze their interconnections and some asymptotic and geometrical properties of these classes.
\end{abstract}
\section{Introduction}

There are  different ways to describe a permutation $\alpha: A\longrightarrow A$ of a set $A$ with $n$ elements, see for example the section "Permutation Statistics" in Stanley's book \cite{stanly}. First we "coordinatise" the set $A$ fixing  a bijection $\beta: A\longrightarrow [n]=\{1,2,\ldots,n\}$ and replace any $A$ permutation $\alpha$ by the $n$-permutation $\pi=\beta\alpha\beta^{-1}$.

The standard representation of such a permutation is the sequence $[\pi]=[\pi(1),\pi(2),$ $\ldots,\pi(n)]$. The first definition of simple permutation is (see \cite{Atk}, \cite{survey}, and Definition 2.3):
\begin{definition}$\pi$ is a \emph{segment-simple} permutation (\emph{s-simple}) if there is no proper connected  subset (or segment) $I\subset [n]$ with a connected image $\pi(I)$.
\end{definition}
Another standard representation of a permutation $\pi$ is given by  a product of  disjoint cycles $(\pi)=(i,\pi(i),\pi^{2}(i)\ldots)\ldots(k,\pi(k),\pi^{2}(k),\ldots)$ (the fixed points are omitted). The \emph{standard representation} (see \cite{stanly}) gives a unique form of a cycle decomposition, $\widehat{\pi}$ : start to write a cycle (the length one cycles are included now) with the biggest element  $k_{i}$ and put the cycles in increasing order of their biggest elements. Now we obtain another sequence of elements in  $[n]:\, \widehat{\pi}=(k_{1}\pi(k_{1})\ldots k_{2}\pi(k_{2})\ldots )$.
\begin{definition}$\pi$ is a \emph{cycle-simple} permutation (\emph{c-simple}) if its cycle decomposition contains at most one cycle of length $\geq 2$.
\end{definition}Now we consider the permutation $\pi$ as an element in the symmetric group $\Sigma_{n}$.

\begin{definition}$\pi$ is a \emph{group-simple} permutation (\emph{g-simple}) if the subgroup generated by $\pi$ is a simple group.
\end{definition}
In the classical Coxeter presentation \cite{Bjorn and Brents} of $\Sigma_{n}$ the generators are the transpositions $\tau_{i}=(i,i+1)$, $i=1,2,\ldots,n-1$. Any element $\pi\in \Sigma_{n}$ can be represented in a unique way as a product of $\tau_{i}$ choosing the smallest such product in the length-lexicographic order given by $\tau_{1}<\tau_{2}<\ldots <\tau_{n-1}$:
$$\pi=(\tau_{k_{1}}\tau_{k_{1}-1}\ldots \tau_{j_{1}})(\tau_{k_{2}}\tau_{k_{2}-1}\ldots \tau_{j_{2}})\ldots(\tau_{k_{s}}\tau_{k_{s}-1}\ldots \tau_{j_{s}})$$ with $1\leq k_{1}<k_{2}<\ldots <k_{s}\leq n-1$ (see \cite{Bjorn and Brents}).

\begin{definition}\label{b-simple def}$\pi$ is a \emph{braid-simple} permutation (\emph{b-simple}) if,  in the  word $\tau_{i_{1}}\tau_{i_{2}}\ldots \tau_{i_{k}}$ representing $\pi$, where $k$ is the length of $\pi$, a Coxeter generator $\tau_{i}$ appears at most once.\end{definition}

There are also graphical descriptions of permutation, see \cite{stanly}. If $(k_{1}\ldots k_{2}\ldots k_{s}\ldots)$ is the standard representation  of $\pi$, then the associated ordered tree $T(\pi)$ is defined inductively as follows: choose as the root of the tree the smallest element $s$ of the sequence, next put on the left the tree of the subsequence  in front of $s$ and on the right  the tree corresponding to the subsequence   behind $s$.
The vertices of $T(\pi)$ are marked from 1 to $n$, 1 is the root of  $T(\pi)$, and along the branches the marking is increasing.
\begin{definition}$\pi$ is a \emph{tree-simple} permutation  (\emph{t-simple}) if all the vertices of the associated tree $T(\pi)$ have degree 1 or 2.
\end{definition}
For instance, the permutation $\pi=\left(
                                     \begin{array}{cccccc}
                                       1 & 2 & 3 & 4 & 5 & 6 \\
                                       4 & 1 & 6 & 2 & 5 & 3 \\
                                     \end{array}
                                   \right)
$ has the following descriptions:

\begin{center}
\begin{picture}(250,100)
\thicklines
\put(-10,0){$\pi=(\tau_{3}\tau_{2}\tau_{1})(\tau_{4})(\tau_{5}\tau_{4}\tau_{3})$}
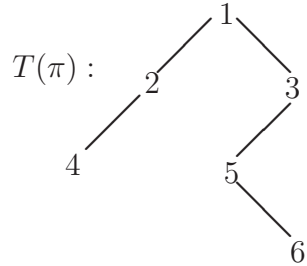
\put(-10,30){$\widehat{\pi}=(4\,2\,1\,5\,6\,3)$}   \put(-10,60){$(\pi)=(4\,2\,1)(6\,3)$}
\put(-10,90){$[\pi]=[4\,1\,6\,2\,5\,3]$}           \put(240,0){$6$}\put(240,10){\line(-1,1){20}}
\put(215,30){$5$}\put(220,40){\line(1,1){20}}      \put(238,62){$3$}\put(240,72){\line(-1,1){20}}
\put(213,90){$1$}\put(210,92){\line(-1,-1){20}}    \put(185,64){$2$}\put(183,63){\line(-1,-1){20}}
\put(155,33){$4$}                                  \put(135,70){$T(\pi):$}
\end{picture}
\end{center}

\medskip

\noindent and this permutation is s- and t-simple, but not b-, c-, or g-simple.

A  last couple of definitions:  we consider a family of subsets $A_{n}\subset\Sigma_{n}$, $n\in\mathbb{N}$.

\begin{definition}The family $(A_{n})$ has \emph{exponential growth} if there are constants $a$, $b>1$ such that card$(A_{n})\geq a\cdot b^{n}$ for any $n\in\mathbb{N}$.
\end{definition}

\begin{definition}The family $(A_{n})$ is \emph{rare} if $\lim\limits_{n\rightarrow\infty}\dfrac{\mbox{card}(A_{n})}{\mbox{card}(\Sigma_{n})}=0.$
\end{definition}

In the next sections we give different characterizations of the corresponding simple subsets of $\Sigma_{n}$: $sS_{n}$, $cS_{n}$, $gS_{n}$, $bS_{n}$, $tS_{n}$ (for instance, in terms of cyclic decomposition), we study enumerative combinatorics of these sets  and growth type  of these subsets. We prove that
\begin{theorem}\label{exp growth}All five simple families $sS_{n}$, $cS_{n}$, $gS_{n}$, $bS_{n}$ and $tS_{n}$ have an exponential growth.
\end{theorem}
\begin{theorem}\label{rare}The four simple families $cS_{n}$, $gS_{n}$, $bS_{n}$ and $tS_{n}$ are rare.
\end{theorem}

We analyze the subgraphs $\Gamma(bS_{n}),\Gamma(cS_{n})$, and $\Gamma(gS_{n})$ of the Cayley graph of $\Sigma_{n}$ (with Coxeter generators), corresponding to these various simple subsets and also the subcomplex $P(bS_{n})$ of the permutahedron $P(\Sigma_{n})$, and the simplicial complexes $B(bS_{n})$, $W(bS_{n})$ associated to Bruhat order and weak order respectively.

The main results are:

\begin{theorem}\label{bSn is conectd}a) $\Gamma(bS_{n})$ is a connected graph which is planar if and only if $n\leq 5$.

\noindent b) $\Gamma(cS_{n})$ is connected if and only if $n\leq 4$.

\noindent c) $\Gamma(gS_{n})$ is connected if and only if $n\leq 3$.
\end{theorem}
Section 5 contains a complete description of connected components of $\Gamma(cS_{n})$ and $\Gamma(gS_{n})$.
\begin{theorem}\label{bSn is contractible}$P(bS_{n})$, $B(bS_{n})$ and $W(bS_{n})$ are  contractible spaces.
\end{theorem}
Our starting point was to understand the permutations corresponding to the simple braids, see~\cite{Barbu-Rehana:02},~\cite{Barbu:Rehana3} for this related notion. As we found in the literature
different notions of simple permutations, the aim of this paper is to characterize in an algebraic and combinatorial way braid-simple permutations and to have a complete picture of the relations between $bS_{n}$ and other families of simple permutations. There are a few natural questions about the other families and these could be interesting problems  in the combinatorics of $\Sigma_{n}$.

\section{Braid-Simple  Permutations}

The set of braid-simple permutations, $bS_{n}\subset\Sigma_{n}$, is the image of the set of simple braids $S_{n}\subset\mathcal{B}_{n}$ through the canonical map $\pi:\,\mathcal{B}_{n}\longrightarrow\Sigma_{n}$ (see~\cite{Barbu-Rehana:02}). Simple braids can be defined as in Definition \ref{b-simple def}, replacing Coxeter  generators by classical Artin generators, therefore the simple braids are square free positive braids or divisors of Garside braid $\Delta_{n}$ (see~\cite{BS},~\cite{Garside:69}): $S_{n}\subset\Div(\Delta_{n})\subset\mathcal{B}_{n}$. As the restriction $\pi:\,\Div(\Delta_{n})\longrightarrow \Sigma_{n}$ is a bijection, we obtain the diagram
\begin{center}
\begin{picture}(150,70)
\thicklines

\put(10,5){$bS_{n}$}   \put(40,3){\LARGE$\hookrightarrow$}  \put(15,50){\vector(0,-1){30}}\put(10,55){$S_{n}$}
\put(40,55){\LARGE$\hookrightarrow$}   \put(70,55){Div$(\Delta_{n})$}
\put(122,55){\LARGE$\hookrightarrow$}\put(155,55){$\mathcal{B}_{n}$}
\multiput(90,50)(70,0){2}{\vector(0,-1){30}}       \put(160,46){\vector(0,-1){30}}
\put(85,5){$\Sigma_{n}$}\put(155,5){$\Sigma_{n}$}  \put(120,5){$=$}
\put(5,30){$\pi$}\put(20,30){$\approx$}
\put(80,30){$\pi$}\put(95,30){$\approx$}
\put(150,30){$\pi$}
\end{picture}
\end{center}

\noindent and therefore  we can use all the results from \cite{Barbu-Rehana:02} and \cite{Barbu:Rehana3}.

We denote by $\sigma_{n,i}$ the number of b-simple permutations of length $i$ in $\Sigma_{n}$ and by $(F_{n})_{n\geq 0}$ the Fibonacci sequence $0,1,1,2,3,5,\ldots$ :

\begin{proposition}(\cite{Barbu:Rehana3}) a)  The numbers $\sigma_{n,i}$ satisfies the recurrences:

1) $\sigma_{1,0}=1$ and $\sigma_{1,i}=0$ for $i\neq 0$;

 2)   $\sigma_{n,i}=\sigma_{n-1,i}+\sigma_{n-1,i-1}+\sigma_{n-2,i-2}+\ldots+\sigma_{n-i,0}$;

 3) $\sigma_{n,i}=2\sigma_{n-1,i-1}+\sigma_{n-1,i}-\sigma_{n-2,i-1}$.

\noindent b)  The cardinality of the set of braid-simple permutations is given by
$$|bS_{n}|=\sigma_{n,0}+\sigma_{n,1}+\ldots+\sigma_{n,n-1}=F_{2n-1}\, .$$
\end{proposition}
We obtain an asymmetric Pascal triangle for $(\sigma_{n,i})$:
\begin{tabbing}
  0 \qquad\qquad\qquad\qquad\qquad\qquad\= 1 \= 2 \= 3 \= 4 \= 5 \= 6 \= 7 \= 8 \=9 \kill
  \> \>  \>  \>  \> 1 \>  \>  \>  \>  \\
  \> \>  \>  \> 1 \>  \> 1 \>  \>  \>  \\
  \> \>  \> 1 \>  \> 2 \>  \> 2 \>  \>  \\
  \> \> 1 \>  \> 3 \>  \> 5 \>  \> 4 \>  \\
 \> 1 \>  \> 4 \>  \> 9 \>  \> 12 \>  \> 8\end{tabbing}
$\indent\indent\indent\indent\indent\indent\indent\indent\indent\indent\indent\indent\ldots\ldots\ldots\ldots\ldots\ldots$

\noindent with $2^{n-2}$ on the last position (for $n\geq 2$).
\begin {cor}The family of braid-simple permutations $bS_{n}$ is rare and has exponential growth.
\end{cor}
\begin{proof}Elementary computations show that
$$F_{2n-1}=c_{1}\Big(\frac{1+\sqrt{5}}{2}\Big)^{2n-1}+c_{2}\Big(\frac{1-\sqrt{5}}{2}\Big)^{2n-1}=
k_{1}\Big(\frac{3+\sqrt{5}}{2}\Big)^{n}+k_{2}\Big(\frac{3-\sqrt{5}}{2}\Big)^{n}.$$
\end{proof}
Now we start to characterize b-simple permutations in terms of cyclic decomposition. First two definitions:
\begin{definition} A subset $I\subset[n]=\{1,2,\ldots,n\}$ is \emph{connected} (or $I$ is a \emph{segment}) if $I=\{i,i+1,\ldots,j\}$ for some $1\leq i\leq j\leq n$. $I$ is \emph{proper} if its cardinality   is not 1 or $n$. A permutation $\pi\in\Sigma_{n}$ is \emph{connected} if all its orbits are connected. For instance $(5\,1\,3\,4\,2)(7\,6)$ is connected but $(4\,2)(5\,1\,3)(7\,6)$ is not.
\end{definition}
\begin{definition} A cycle $(k_{1}\,k_{2}\,\ldots k_{s})$ is \emph{unimodal} if there is an index $1\leq m\leq s$ such that $k_{1}>k_{2}>\ldots >k_{m}<k_{m+1}<\ldots <k_{s}$. A permutation $\pi$ is \emph{unimodal} if all its cycles are unimodal.
 \end{definition}
 In this definition we use the standard representation convention, $k_{1}$ is the largest element of the cycle. If we start to write the cycle with the smallest element, we find the sequence $(k_{m},k_{m+1},\ldots,k_{s},k_{1},k_{2},\ldots,k_{m-1})$, which is unimodal in the usual sense (first increasing, next decreasing). In fact we need "cyclic unimodal"  sequences, hence the usual definition and our definition coincide. We will use the notations
 $D(k,j)=\tau_k\tau_{k-1}\ldots \tau_{j+1}\tau_j $ for $1\leq j\leq k\leq n-1$ and also $ D(K_*,J_*)= D(k_1,j_1)D(k_2,j_2)...D(k_s,j_s)$, where
 $1\leq k_1<k_2<\ldots <k_s\leq n-1$, $ j_{a}\leq k_{a}$, and $ 1\leq s \leq n-1 $. Also $D(k)$ is a short notation
for $D(k,k)=\tau_k$.

 The new result of this section is:
 \begin{proposition}A permutation $\pi\in\Sigma_{n}$ is braid-simple if and only if $\pi$ is connected and unimodal.
 \end{proposition}
 \begin{proof}A product $D(k,j)=\tau_{k}\tau_{k-1}\ldots \tau_{j+1}\tau_{j}$ has the cycle representation $(k+1,k,\ldots,$
 $j+1,j)$. A b-simple permutation has the canonical representation as a product of Coxeter generators
 $$D(K_*,J_*)= D(k_1,j_1)D(k_2,j_2)\ldots D(k_s,j_s) $$
 where $1\leq j_{1}\leq k_{1}< j_{2}\leq k_{2}<\ldots <j_{s-1}\leq k_{s-1}< j_{s}\leq k_{s}\leq n-1$. This permutation is a cycle if and only if $j_{2}=k_{1}+1, \ldots, j_{s}=k_{s-1}+1$ and in this case its cycle representation is
 $$(k_{s}+1,k_{s},\ldots,j_{s}+1,k_{s-1},k_{s-1}-1,\ldots,j_{s-1}+1,\ldots ,k_{1},k_{1}-1,\ldots,j_{1},j_{2},\ldots,j_{s})$$
 (if one of the factors $D(k_{i},j_{i})$ contains only one generator, $\tau_{k_{i}}=\tau_{j_{i}}$, then this generator appears only in the increasing part of this connected unimodal cycle).

 Conversely, a unimodal connected cycle can be written as a product of $D's$ factors starting with the leftmost segment containing the minimal element in the cycle, next reading  the segments (from right to left), adding at their ends  one element from the increasing part of the unimodal sequence and writing these augmented segments from left to right. As an example, the connected unimodal cycle $(13,12,9,8,7,5,3,2,1,4,6,10,11)$ corresponds to the product of Coxeter generators $D(3,1)D(5,4)D(9,6)D(10)D(12,11)$. If the product $D(K_{*},J_{*})$ contains $c-1$ "jumps"
 of the form $k_{i}+1<j_{i+1}$, then its cycle decomposition has $c$ disjoint cycles, and all of them are connected and unimodal.\end{proof}
 Direct consequences of the proof and of Proposition \ref{g-simple characterization} are the next two characterizations of the braid-simple permutations which are cycle-simple and group-simple, respectively.
 \begin{cor}The following properties of a permutation $\pi\in \Sigma_{n}$ are equivalent:

 i) $\pi\in bS_{n}\bigcap cS_{n}$;

 ii) $\pi=Id$ or the cycle representation of $\pi$ contains a unique cycle and this is connected and unimodal;

 iii)  $\pi=Id$ or the Coxeter representation of $\pi$ is $D(K_{*},J_{*})$ where $k_{i}+1=j_{i+1}$ for any $i$.
  \end{cor}
  \begin{cor}The following properties of a permutation $\pi\in \Sigma_{n}$ are equivalent:

 i) $\pi\in bS_{n}\bigcap gS_{n}$;

 ii) $\pi=Id$ or there exist a prime $p$ and $q$ ($pq\leq n$) such that the cycle decomposition of $\pi$ contains $q$ cycles of length $p$, and  all of them are connected and unimodal;

 iii) $\pi=Id$ or there exist a prime $p$ and $q$ ($pq\leq n$) such that the Coxeter representation of $\pi$ is $D(K_{*},J_{*})=D(K_{*}^{1},J_{*}^{1})D(K_{*}^{2},J_{*}^{2})\ldots D(K_{*}^{a},J_{*}^{a}),$ where $D(K_{*}^{a},J_{*}^{a})=D(k_{1}^{a},j_{1}^{a})\ldots D(k_{s_{a}}^{a},j_{s_{a}}^{a})$ with $k^{a}_{i}+1=j^{a}_{i+1}$, $k^{a}_{s_{a}}-j^{a}_{1}+2=p$, and $k^{a}_{s_{a}}+1 < j^{a+1}_{1}$.
\end{cor}

 \section{cyclic-simple and group-simple permutations}
 Let start with the remark that cyclic-simple and group-simple notions are invariant under conjugation, i.e. these notions does not depend on the "coordinatization" of a finite set $A$ with $n$ elements. All the other classes, b-, s-, and t-simple permutations, depend on the coordinatization of the set $A$.

  The next characterization is obvious:
 \begin{proposition}\label{g-simple characterization} A permutation $\pi\in\Sigma_{n}$ is g-simple if and only if $\pi=$identity or there is a prime number $p$ and a positive integer $k$ ($pk\leq n$) such that $\pi$  has $k$ cycles of length $p$ (and other elements are fixed).  \end{proposition}
 Enumerative combinatorics of $cS_{n}$ and $gS_{n}$ is simple:
 \begin{proposition}a) The number of permutations $\pi\in\Sigma_{n}$, product of $k$ cycles of  equal length $l$, is given by:
 $$\frac{n!}{k!(n-kl)!l^{k}};$$
 b) the set of cyclic-simple permutations $cS_{n}$ has cardinality:
 $$1+\sum\limits_{l=2}^{n}\frac{n!}{l\cdot(n-l)!};$$
 c) the set of group-simple permutations $gS_{n}$ has cardinality:
 $$1+\sum\limits_{p\,\,prime}\,\,\sum\limits_{k=1}^{\lfloor\frac{n}{p}\rfloor}\frac{n!}{k!(n-kp)!p^{k}}.$$
 \end{proposition}
 \begin{proof}a) There are $\left(
   \begin{array}{c}
       n \\
       l\,l\,\ldots l\,\,n-kl \\
       \end{array}
       \right) $ (with $l$ repeated $k$ times) choices for the sequence of $k$ orbits  of length  $l$,
 $\dfrac{1}{k!}
 \left(    \begin{array}{c}
        n \\
        l\,l\,\ldots l\,\,n-kl \\
       \end{array}
       \right) $
 choices for the set of these $k$-orbits; for each orbit of length $l$, there are $(l-1)!$
 cycles, hence the result.

 b) and c) are consequences of a).  \end{proof}
 \begin{cor}The two families $(cS_{n})$ and $(gS_{n})$ are rare with exponential growth.
 \end{cor}
\begin{proof}In the sum representing $|cS_{n}|$ the last term is $(n-1)!\sim \sqrt{2\pi (n-1)}(\frac{n-1}{e})^{n-1}$. By Bertrand postulate \cite{Bertrend}, there is a prime $p$ between $\lfloor\frac{n}{2}\rfloor$ and $n$, therefore the last term in the double sum representing $|gS_{n}|$, $\dfrac{n!}{(n-p)! p}$, is greater than $\dfrac{n!}{\lfloor\frac{n}{2}\rfloor! n} \sim \dfrac{\sqrt{2}}{n}\Big(\dfrac{2n}{e}\Big)^{\lfloor\frac{n}{2}\rfloor}$.

$$\text{ We split the sum }\frac{1}{n!}+\sum\limits_{l=2}^{n}\frac{1}{(n-1)!l}=\frac{|cS_{n}|}{|\Sigma_{n}|}\text{ in two parts}\indent\indent\indent\indent\indent\indent\indent\indent$$

$$C_{1}=\frac{1}{n!}+\sum\limits_{l=2}^{\lfloor\frac{n}{2}\rfloor}\frac{1}{(n-l)!l}\leq \frac{\lfloor\frac{n}{2}\rfloor}{2 \lfloor\frac{n}{2}\rfloor !}$$
 and
$$C_{2}=\sum\limits_{l=\lfloor\frac{n}{2}\rfloor +1}^{n}\frac{1}{(n-l)!l}\leq \dfrac{1}{\lfloor\frac{n}{2}\rfloor} (1+\frac{1}{1!}+\frac{1}{2!}+\ldots+\frac{1}{n!})<\frac{e}{\lfloor\frac{n}{2}\rfloor}$$
hence $\lim \dfrac{|cS_{n}|}{n!}=0.$

In a similar way we can evaluate the quotient $\dfrac{|gS_{n}|}{|\Sigma_{n}|}$: for primes $p$ in the 
interval $(\lfloor\frac{n}{2}\rfloor,n]$, we have
$$G_{2}=\sum\limits_{\lfloor\frac{n}{2}\rfloor< p\leq n}\sum\limits_{k=1 }^{\lfloor\frac{n}{p}\rfloor}\frac{1}{k! (n-kp)! p^{k}}
=\sum\limits_{\lfloor\frac{n}{2}\rfloor < p \leq n}\frac{1}{(n-p)! p}
<\frac{1}{\lfloor\frac{n}{2}\rfloor}(1+\frac{1}{1!}+
\frac{1}{2!}+\ldots+\frac{1}{n!})<\frac{e}{\lfloor\frac{n}{2}\rfloor}
,$$
otherwise $p\leq \lfloor\frac{n}{2}\rfloor$ and we split in two parts the contribution of $p$:
$$g_{p}=\sum\limits_{k=1}^{\lfloor\frac{n}{p}\rfloor}\frac{1}{k! (n-kp)!p^{k}}=\sum\limits_{k\leq  \lfloor\frac{n}{2p}\rfloor}\frac{1}{k! (n-kp)!p^{k}}+\sum\limits_{k > \lfloor\frac{n}{2p}\rfloor}\frac{1}{k! (n-kp)!p^{k}}=g'_{p}+g''_{p}.$$

For the first sum we find
$$g'_{p}=\sum\limits _{k\leq \lfloor\frac{n}{2p}\rfloor}\frac{1}{k!(n-kp)!p^{k}}\leq \lfloor\frac{n}{2p}\rfloor
\frac{1}{1!\lfloor\frac{n}{2}\rfloor ! 2}
\leq \frac{\lfloor\frac{n}{4}\rfloor}{2\lfloor\frac{n}{2}\rfloor !}
$$
and for the second sum we find
$$g''_{p}=\sum\limits _{k >\lfloor\frac{n}{2p}\rfloor}\frac{1}{k!(n-kp)!p^{k}}\leq \lfloor\frac{n}{2p}\rfloor \frac{1}{\lfloor\frac{n}{2p}\rfloor ! p^{\lfloor\frac{n}{2p}\rfloor}}<\frac{constant}{n^{2}}$$
(for the last inequality we use $\dfrac{n^{2}}{c}< \left(
                                                           \begin{array}{c}
                                                             \lfloor\frac{n}{2p}\rfloor \\
                                                             2\\
                                                           \end{array}
                                                         \right)
$ $(p-1)^{2}< p^{\lfloor\frac{n}{2p}\rfloor}$, with $c$ a constant in $n$ and $p$). Finally the inequality for $G_{2}$ and
$$G_{1}=\frac{1}{n!}+\sum\limits _{p \leq \lfloor\frac{n}{2}\rfloor}(g'_{p}+g''_{p})< \frac{1}{n!}+\lfloor\frac{n}{2}\rfloor (\frac{\lfloor\frac{n}{4}\rfloor}{2\lfloor\frac{n}{2}\rfloor !}+\frac{constant}{n^{2}})$$
shows that $\lim\dfrac{|gS_{n}|}{n!}=0$.
 \end{proof}

\section{Segment-simple and tree-simple permutations}
In this section we finish the proofs of Theorems  \ref{exp growth}  and \ref{rare} and we characterize the permutations from the intersections $bS_{n}\bigcap sS_{n}$
and $bS_{n}\bigcap tS_{n}$.

The asymptotics of segment-simple permutations is given by the next result:
\begin{theorem}(\cite{Atk}) The cardinality of the set of segment-simple permutations is
$$|sS_{n}|=\frac{n!}{e^{2}}\Big(1-\frac{4}{n}+\frac{2}{n(n-1)}+O(n^{-3})\Big).$$
\end{theorem}
Enumerative combinatorics  of tree-simple permutations is elementary:

\begin{proposition}The cardinality of the set of tree-simple permutations is,  for $n\geq 2$,
$$|tS_{n}|=2^{n-2}+4^{n-2}. $$
\end{proposition}
\begin{proof}For $n\geq 2$, first we count the number of simple oriented rooted marked trees: marks are from 1 to $n$, the number 1 is the root, at this point we could have a left  branch, or a right branch, or both, at any other node we have at most one branch (a left one or a right one), and all the marks on the branches are in increasing order. If there is only one branch at the root 1, we have $2^{n-1}$ choices for left-right orientations at the nodes $1,2,\ldots, n-1$. If there are $k$ nodes on the left branch and $(n-k-1)$ nodes on the right branch ($1\leq k \leq n-2$), we have $\left(                                                                                                \begin{array}{c}                                                                                                                            n-1 \\
k \\
\end{array}                                                                                                                                 \right)  $
ways to choose the increasing marking on the left branch, $2^{k-1}$ possibilities to choose orientations at the first $k-1$ nodes of this branch, and $2^{n-k-2}$ possibilities to choose left-right orientations at the first $n-k-2$ nodes of the right branch. Therefore the total number is ($n\geq 2$):
$$2^{n-1}+\sum\limits _{k=1}^{n-2}\left(
                                    \begin{array}{c}
                                      n-1 \\
                                      k \\
                                    \end{array}
                                  \right)
2^{k-1}\cdot 2^{n-k-2}=2^{n-1}+(2^{n-1}-2)2^{n-3}=2^{n-2}+4^{n-2}.$$
For $n=1$, the last formula (with $\sum=0$) gives the result $|tS_{1}|=1.$
\end{proof}
In the next table one can see few values of the numbers of $*$-simple permutations. For large values of $n$,
there are more segment-simple permutations than in the union of the other four classes.
 \begin{center}
 \begin{tabular}{|c||c|c|c|c|c|c|}
   \hline
   $n$ & $|bS_{n}|$ & $|cS_{n}|$ & $|gS_{n}|$ & $|sS_{n}|$ & $|tS_{n}|$ & $|\Sigma_n|$ \\\hline
   1 & 1 & 1 & 1 & 1 & 1 & 1 \\\hline
   2 & 2 & 2 & 2 & 2 & 2 & 2 \\\hline
   3 & 5 & 6 & 6 & 0 & 6 & 6 \\\hline
   4 & 13 & 21 & 18 & 2 & 20 & 24 \\\hline
   5 & 34 & 85 & 70 & 6 & 72 & 120 \\\hline
   6 & 89 & 410 & 300 & 46 & 272 & 720 \\\hline

 \end{tabular}
 \end{center}

 \medskip
A cycle \emph{contains three consecutive elements} if it contains a subsequence $(\ldots,i,i+1,i+2,\ldots)$ or $(\ldots,i+2,i+1,i,\ldots)$. For example, the cycle $(7431256)$ contains three consecutive elements $(567)$, but $(312)$ are not consecutive.
\begin{proposition}The following properties of a permutation $\pi\in\Sigma_{n}$ are equivalent:

i) $\pi\in bS_{n}\bigcap sS_{n}$;

ii) $\pi$ is an unimodal cycle of length $n$ without  three consecutive elements.
\end{proposition}
\begin{proof}If $\pi$ is braid-simple in $\Sigma_{n}$ it is connected and unimodal, therefore if $\pi$ has a cycle of length $2\leq q \leq n-1$, a segment of length $q$ is $\pi$ invariant. This shows that $\pi \in bS_{n}\bigcap sS_{n}$ should be an $n$-cycle. If there are three consecutive elements, $(i-1,i,i+1)$, $(i+1,i,i-1)$, then $\pi$(segment of length 2)  is also a segment: $\pi[i,i\pm1]=[i,i\mp 1]$.

Conversely, if $\pi$ is a unimodal $n$-cycle without three consecutive elements, then we have to show that $\pi$ is segment-simple. Suppose that there are two segments $a\leq b$, $c\leq d$, $b-a=d-c=L, 2\leq L\leq n-1$, such that $\pi[a,b]=[c,d]$. The segment $[a,b]$ can not be contained entirely in the decreasing part $(n\ldots 1)$ or in the increasing part $(1\ldots k)$ of the unimodal cycle $\pi$ (if $L\geq 3$, then we have three consecutive elements; if $L=2$, then $a\rightarrow a+1$, and now $\pi$[a,a+1]=[a+1,a+2], therefore again we find three consecutive elements, and similarly, if $a+1\rightarrow a$ we find $a\rightarrow a-1$, contradiction). Let us denote by $d_{M}$ and $d_{m}$ the maximal and minimal elements in the decreasing part of the cycle $\pi$ contained in the interval $[a,b]$, and by $i_{m}$ and $i_{M}$ the minimal and the maximal elements of $[a,b]$ contained in the increasing part of $\pi$:
$$\pi=(n\ldots d_{M}\ldots d_{m}\ldots 1\ldots i_{m}\ldots i_{M}\ldots k).$$
The value $\pi(d_{m})$ is less than $d_{m}$ and not in $[a,b]$, therefore $\pi(d_{m})< a$ and $\pi(i_{M})$ is greater than $i_{M}$ and not in $[a,b]$, therefore $\pi(i_{M})> b$, but these inequalities give a contradiction
$$\pi(i_{M})-\pi(d_{m})>b-a=L=c-d \geq \pi(i_{M})-\pi(d_{m}).$$
The special cases $d_{m}=1$, or $i_{m}=1$, or $d_{m}=i_{m}=1$ can be settled in a similar way.
\end{proof}
To represent the rooted tree of a permutation, $T(\pi)$, without self intersections, we make the following "geometrical" conventions:

a) all the edges make  an angle equal to $\dfrac{\pi}{4}$  with the vertical direction;

b) the edges starting at the root 1 have the length 1, the next edges (at most four) have length $\dfrac{1}{2}$, the next (at most eight) edges have length $\dfrac{1}{4}$, and so on.

In the tree $T(\pi)$ the ending points of the left and right branches are denoted by $M_{L}$ and $M_{R}$ (they are the maximal marks on their branches).

\begin{proposition}The following properties of a permutation $\pi\in\Sigma_{n}$ are equivalent:

i) $\pi\in bS_{n}\bigcap tS_{n}$;

ii)  there is at most one right angle (different from 1) and this could be only on the right branch. In the case of a right angle on the right branch at the point $M$, we must have $M_{L}< M$.
\end{proposition}
\begin{proof} If $\pi$ is braid-simple, then $\pi$ is connected unimodal, in particular the cycle containing 1 contains only the numbers from 1 to $a$: ($a=a_{1}, a_{2} \ldots ,1,\ldots,a_{k}$), therefore the first part of the tree $T({\pi})$ has (at most) a right angle at 1. Because $\pi$ 

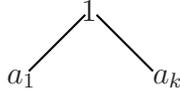
\begin{figure}[hbt]
\begin{center}
\begin{pspicture}(0,0.8)(2,0)
\psline(0.75,0.75)(0,0)\rput(0.8,0.8){$1$}
\psline(0.9,0.75)(1.65,0)\rput(-0.1,-0.1){$a_{1}$}\rput(1.85,-0.1){$a_{k}$}
\end{pspicture}
\end{center}
\end{figure}

\noindent is tree-simple, the next value, $a_{1}+1$, should be either the successor of $a_{k}$ (and hence a fixed point of $\pi$) or the last one in the standard representation; in the last case, the cycle containing $a_{1}+1$ contains also $n$ (connectedness condition) and must be $(n,n-1,\ldots,a_{1}+1)$ ($\pi$ is unimodal). Hence we have a unique right angle  at $a_{1}+1$. If $\pi(a_{1}+1)=a_{1}+1$, the same argument says that the structure of $\pi$ is:
$$\pi=(a_{1}\ldots 1 \ldots a_{k})(a_{1}+1)\ldots(a_{1}+i-1)(n\,n-1 \ldots a_{1}+i)$$
or
$$\pi=(a_{1}\ldots 1 \ldots a_{k})(a_{1}+1)(a_{1}+2)\ldots n,$$
and the corresponding trees have one or none right angles on the right branch and they have no right angles on the left branch. Using the same argument, in the  particular case of $\pi(1)=1$, $T(\pi)$ has only the right branch with at most one right angle. If $a_{1}=n$, then $\pi=(n,n-1,\ldots,2,1)$ has only a straight left branch.
\begin{figure}[hbt]
\begin{pspicture}(0,-1.5)(6,1.5)
\rput(-1.5,0){\rput(0.5,0.5){$M_L=a_{1}$}\psline(0.6,0.6)(1.3,1.3)
\rput(1.4,1.4){$1$}\psline(1.5,1.3)(2.2,0.6)\rput(2.5,0.5){$a_{1}+1$}
\psline(2.3,0.4)(3,-0.3)\rput(3.4,-0.4){$n=M_R$}
\rput(-1,1){$T_{1}:$}}
\rput(5,0){\rput(0.5,0.5){$M_L=a_{1}$}\psline(0.6,0.6)(1.3,1.3)
\rput(1.4,1.4){$1$}\psline(1.5,1.3)(2.2,0.6)\rput(2.5,0.5){$a_{1}+1$}
\psline(2.3,0.4)(3,-0.3)\rput(3.7,-0.4){$a_{j}=M$}
\psline(3.1,-0.5)(2.4,-1.2)\rput(2.7,-1.3){$n=M_R$}\rput(-1,1){$T_{2}:$}}
\end{pspicture}
\end{figure}

Conversely, the  tree $T_{1}$  corresponds to 
$$\pi_{1}=(a_{1}\ldots 1 \ldots a_{k})(a_{1}+1)\ldots (n),$$
and the tree $T_{2}$ corresponds to the simple permutation ($a_j$ is greater than $a_1$)
$$\pi_{2}=(a_{1}\ldots 1 \ldots a_{k})(a_{1}+1)\ldots (a_{j}-1)(n,\,n-1 \ldots a_{j}),$$
and both are b- and t-simple; we have  similar results for trees without right angles.
\end{proof}
\begin{cor}The intersection $S_{n}=bS_{n}\bigcap cS_{n}\bigcap gS_{n} \bigcap sS_{n}\bigcap tS_{n}$ is non empty if and only if  $n\leq 2$ or $n$ is a prime number greater or equal to 5.
\end{cor}
\begin{proof}For $n=1$ and $n=2$, the intersection $S_{n}=bS_{n}\bigcap cS_{n}\bigcap gS_{n} \bigcap sS_{n}\bigcap tS_{n}$ is the entire set $\Sigma_{n}$. For $n=3$, $sS_{3}=\emptyset$ and for $n=4$, see the first picture in section 5. Because $\pi$ is braid-simple and segment-simple, $\pi$ should be an $n$-cycle; $\pi$ is group-simple implies that $n$ should be a prime number $\geq 5$.

Conversely, if $p$ is prime and  $\geq 5$, the permutation $$\pi=(p,p-2,p-4,\ldots,5,3,1,2,4,\ldots,p-3,p-1)$$ is 
an element of the intersection $S_{p}$.
\end{proof}

\section{Geometry of Simple Permutations}

If $A$ is a subset of $\Sigma_{n}$, the corresponding \emph{subgraph} $\Gamma(A)$ of the Cayley graph $\Gamma(\Sigma_{n})$ has $A$ as a set of vertices and there is an edge $a$---$b$ between two elements of $A$  if $b=a\tau_{i}$ for some Coxeter generator $\tau_{i}=(i,i+1)$  ($b=a\tau_{i}$ is equivalent to $a=b\tau_{i}$). In a similar way is defined the \emph{subcomplex} $P(A)$ of the permutahedron $P(\Sigma_{n})$: the vertices of $P(A)$  are elements of $A$ and a face of permutahedron is a face of $P(A)$ if and only if all its vertices are in $A$. It is obvious that $\Gamma(A)$ is the 1-dimensional skeleton of $P(A)$. In this last section we analyze subgraphs $\Gamma(*S_{n})$  corresponding to some simple subsets $*S_{n}$ of the symmetric group $\Sigma_{n}$ and the subcomplex $P(bS_{n})$.

Other geometric objects associated  to the Coxeter presentation of $\Sigma_{n}$ are the simplicial complexes corresponding to the natural poset structures on $\Sigma_{n}$: Bruhat order and (right) weak order (see \cite{Bjorn and Brents}). The permutations $\alpha$ and $\beta$ are consecutive in these order relations, $\alpha< \beta$, if and only if length$(\beta)=$length$(\alpha)+1$ and for some transpositions $\beta=\alpha(i,j)$ in the Bruhat order and $\beta=\alpha\tau_{i}$ in the weak order. We denote by $B(\Sigma_{n})$ and $W(\Sigma_{n})$ the corresponding simplicial complexes associated to these two posets and by $B(*S_{n})$ and $W(*S_{n})$ the subcomplexes corresponding to the $*$-simple permutations in $\Sigma_{n}$.

\medskip

\noindent \emph{Proof of Theorem \ref{bSn is conectd}.} a)  The b-simple part $\Gamma(bS_{n})$ is connected by the very definition of a braid simple permutation: $D(K_{*}, J_{*})$ can be joined with the unit element by the geodesic

$\,\,\,\,\,\,\,\,\,\,Id$ ---$\tau_{k_{1}}$---$\tau_{k_{1}}\tau_{k_{1}-1}$---$ \ldots$---$D(k_{1},j_{1})$---$D(k_{1},j_{1})\tau_{k_{2}}$---$ \ldots$---$D(K_{*},J_{*}).$

The planarity of $\Gamma(bS_{6})$ comes from the next figure, and one can find a $K_{3,3}$ subgraph of $\Gamma(bS_{7})$ on the next page.

b) This is a consequence of Proposition \ref{cycle graph is connected}.

c) This is a consequence of Proposition \ref{group graph is connected}.\mbox{}\quad  \mbox{}\quad \mbox{}\quad $\Box$

\begin{center}
\begin{picture}(400,430)
\put(198,188){$\bullet$\tiny$\underline{e}^*$}
\put(218,198){$\bullet$}\put(211,204){\tiny$\underline{1}^*$}  \put(188,168){$\bullet$\tiny$\underline 2$}
\put(198,208){$\bullet$}\put(202,204){\tiny$\underline 3$}     \put(208,168){$\bullet$}\put(202,168){\tiny$\underline 4$}
\put(178,198){$\bullet$}\put(185,203){\tiny$\underline{5}^*$}  \put(238,198){$\bullet$}\put(230,204){\tiny$\underline{12}^*$}
\put(178,148){$\bullet$\tiny$\underline{21}^*$}                \put(218,218){$\bullet$}\put(222,215){\tiny$13$}
\put(158,159){$\bullet$}\put(162,157){\tiny$\underline {23}$}  \put(188,228){$\bullet$\tiny$\underline {32}$}
\put(228,178){$\bullet$}\put(218,178){\tiny$14$}               \put(198,148){$\bullet$\tiny$24$}
\put(208,228){$\bullet$\tiny$\underline {34}$}                 \put(238,159){$\bullet$}\put(230,157){\tiny$\underline {43}$}
\put(198,368){$\bullet$}\put(198,375){\tiny$15^*$}             \put(168,178){$\bullet$\tiny$25$}
\put(178,218){$\bullet$\tiny$35$}                              \put(218,148){$\bullet$\tiny$\underline {45}$}
\put(158,198){$\bullet$}\put(158,205){\tiny$\underline{54}^*$} \put(268,198){$\cdot$}\put(260,205){\tiny$\underline{123}^*$}
\put(258,238){$\cdot$\tiny$\underline{132}^*$}                 \put(158,128){$\cdot$\tiny$\underline{213}^*$}
\put(178,248){$\cdot$\tiny$\underline{321}^*$}                 \put(248,178){$\cdot$}\put(250,185){\tiny$124$}
\put(228,236){$\cdot$\tiny$134$}                               \put(138,128){$\cdot$\tiny$\underline {234}$}
\put(198,248){$\cdot$\tiny$\underline {324}$}                  \put(188,128){$\cdot$\tiny$214$}
\put(308,97){$\cdot$\tiny$143$}                                \put(198,107){$\cdot$\tiny$\underline {243}$}
\put(258,128){$\cdot$\tiny$\underline {432}$}                  \put(278,217){$\cdot$}\put(265,221){\tiny$125^*$}
\put(88,98){$\cdot$}\put(75,98){\tiny$235$}                    \put(118,217){$\cdot$}\put(120,221){\tiny$154^*$}
\put(218,248){$\cdot$\tiny$\underline {345}$}                  \put(128,198){$\cdot$}\put(125,205){\tiny$\underline{543}^*$}
\put(68,58){$\cdot$}\put(65,50){\tiny$215^*$}                  \put(168,238){$\cdot$\tiny$325$}
\put(208,128){$\cdot$\tiny$245$}                               \put(138,238){$\cdot$\tiny$\underline {354}$}
\put(198,348){$\bullet$}\put(203,354){\tiny$135$}              \put(328,57){$\cdot$\tiny$145$}
\put(148,178){$\cdot$}\put(136,184){\tiny$254$}                \put(238,129){$\cdot$\tiny$\underline {435}$}
\put(290,178){$\bullet$\tiny$\underline{1243}^*$}              \put(298,188){$\bullet$}\put(305,190){\tiny$\underline{1234}^*$}
\put(268,258){$\bullet$\tiny$\underline{1324}^*$}              \put(318,78){$\bullet$}\put(327,78){\tiny$\underline{1432}^*$}
\put(139,108){$\bullet$\tiny$\underline{2134}^*$}              \put(188,88){$\bullet$}\put(167,90){\tiny$\underline{2143}^*$}
\put(188,268){$\bullet$}\put(168,268){\tiny$\underline{3214}^*$} \put(278,108){$\bullet$\tiny$\underline{4321}^*$}
\put(308,208){$\cdot$\tiny$1235^*$}                            \put(158,257){$\cdot$\tiny$3215^*$}
\put(118,108){$\bullet$}\put(103,114){\tiny$\underline{2345}^*$} \put(108,178){$\bullet$}\put(90,178){\tiny$\underline{2543}$}
\put(128,258){$\bullet$\tiny$\underline{3254}$}                \put(98,188){$\bullet$}\put(77,190){\tiny$\underline{5432}^*$}
\put(78,78){$\bullet$}\put(58,78){\tiny$\underline{2354}$}     \put(278,158){$\cdot$}\put(276,165){\tiny$1245$}
\put(198,387){$\bullet$}\put(205,390){\tiny$1254$}             \put(198,28){$\bullet$\tiny$2145$}
\put(118,158){$\cdot$}\put(104,164){\tiny$2154^*$}             \put(238,258){$\cdot$\tiny$1345$}
\put(99,258){$\cdot$}\put(102,258){\tiny$1354$}                \put(268,68){$\cdot$}\put(265,63){\tiny$1435$}
\put(88,208){$\cdot$}\put(67,208){\tiny$1543^*$}               \put(208,88){$\bullet$\tiny$\underline{2435}$}
\put(208,268){$\bullet$}\put(210,270){\tiny$\underline{3245}$} \put(258,108){$\bullet$}\put(253,102){\tiny$\underline{4325}$}
\put(298,258){$\cdot$\tiny$1325^*$}                            \put(128,68){$\cdot$\tiny$2135^*$}
\put(328,178){$\cdot$\tiny$\underline{12345}^*$}               \put(198,408){$\cdot$\tiny$\underline{12543}$}
\put(200,330){\circle{6}\tiny$\underline{13254}^*$}\put(197,327){$\star$}
\put(58,197){$\cdot$}\put(32,198){\tiny$\underline{15432}^*$}
\put(91,71){\circle{6}}\put(94,76){\tiny$\underline{21354}^*$}\put(88,68){$\star$}
\put(118,278){$\cdot$\tiny$\underline{32154}^*$}               \put(68,177){$\cdot$}\put(42,178){\tiny$\underline{54321}^*$}
\put(88,157){$\cdot$}\put(62,158){\tiny$\underline{21543}^*$}  \put(308,158){$\cdot$\tiny$\underline{12435}^*$}
\put(278,278){$\cdot$}\put(252,278){\tiny$\underline{13245}^*$}
\put(311,71){\circle{6}}\put(317,68){\tiny$\underline{14325}^*$}\put(308,68){$\star$}
\put(118,88){$\cdot$\tiny$\underline{21345}^*$}                \put(198,288){$\cdot$\tiny$\underline{32145}^*$}
\put(278,88){$\cdot$}\put(269,80){\tiny$\underline{43215}^*$}  \put(201,71){\circle{6}\tiny$\underline{21435}^*$}
\put(198,68){$\star$}                                          \put(338,198){$\cdot$\tiny$\underline{12354}^*$}
\put(180,370){\line(1,0){40}}   \put(180,350){\line(1,0){40}}    \put(100,221){\line(1,0){20}}
\put(160,221){\line(1,0){20}}   \put(220,221){\line(1,0){20}}    \put(280,221){\line(1,0){20}}
\put(130,200){\line(1,0){50}}   \put(220,200){\line(1,0){50}}    \put(110,180){\line(1,0){60}}
\put(230,180){\line(1,0){60}}   \put(90,160){\line(1,0){30}}     \put(280,160){\line(1,0){30}}
\put(90,70){\line(1,0){40}}     \put(270,70){\line(1,0){40}}     \put(70,60){\line(1,0){65}}
\put(265,60){\line(1,0){65}}
\multiput(20,260)(360,0){2}{\line(0,1){40}}   \multiput(40,260)(320,0){2}{\line(0,1){40}}
\multiput(60,260)(280,0){2}{\line(0,1){40}}   \multiput(70,265)(260,0){2}{\line(0,1){30}}
\multiput(80,270)(240,0){2}{\line(0,1){20}}   \multiput(100,260)(200,0){2}{\line(0,1){20}}
\multiput(180,200)(40,0){2}{\line(0,1){20}}   \put(200,190){\line(0,1){20}}
\multiput(160,130)(80,0){2}{\line(0,1){30}}   \multiput(140,110)(120,0){2}{\line(0,1){20}}
\multiput(120,90)(160,0){2}{\line(0,1){20}}   \multiput(190,90)(20,0){2}{\line(0,1){40}}
\multiput(200,30)(0,80){2}{\line(0,1){40}}    \put(200,390){\line(0,1){20}}
\put(200,350){\line(0,1){20}}  
\put(220,250){\line(2,1){60}}   \put(180,250){\line(-2,1){60}}
\put(210,230){\line(2,1){60}}   \put(190,230){\line(-2,1){60}}
\put(200,210){\line(2,1){120}}  \put(200,210){\line(-2,1){120}}
\put(240,220){\line(2,1){90}}   \put(160,220){\line(-2,1){90}}
\put(220,200){\line(2,1){120}}  \put(180,200){\line(-2,1){120}}
\put(240,200){\line(2,1){120}}  \put(160,200){\line(-2,1){120}}
\put(300,220){\line(2,1){80}}   \put(100,220){\line(-2,1){80}}
\put(300,280){\line(-2,1){100}} \put(100,280){\line(2,1){100}}
\put(320,290){\line(-2,1){120}} \put(80,290){\line(2,1){120}}
\put(330,295){\line(-2,1){110}} \put(70,295){\line(2,1){110}}
\put(340,300){\line(-2,1){140}} \put(60,300){\line(2,1){140}}
\put(360,300){\line(-2,1){140}} \put(40,300){\line(2,1){140}}
\put(380,300){\line(-2,1){180}} \put(20,300){\line(2,1){180}}
\put(260,240){\line(1,2){20}}   \put(220,220){\line(1,2){20}}
\put(200,210){\line(1,2){20}}   \put(190,230){\line(1,2){20}}
\put(180,250){\line(1,2){20}}   \put(220,250){\line(-1,2){20}}
\put(210,230){\line(-1,2){20}}  \put(200,210){\line(-1,2){20}}
\put(180,220){\line(-1,2){20}}  \put(140,240){\line(-1,2){20}}
\put(200,190){\line(2,1){20}}   \put(200,190){\line(-2,1){20}}
\put(200,190){\line(-1,-2){65}} \put(200,190){\line(1,-2){65}}
\put(150,180){\line(-3,-2){120}} \put(250,180){\line(3,-2){120}}
\put(170,180){\line(-3,-2){120}} \put(230,180){\line(3,-2){120}}
\put(170,180){\line(-1,-1){80}}  \put(230,180){\line(1,-1){80}}
\put(190,170){\line(-4,-1){30}}  \put(210,170){\line(4,-1){30}}
\put(160,160){\line(-1,-1){60}}  \put(240,160){\line(1,-1){60}}
\put(160,160){\line(-2,-3){20}}  \put(240,160){\line(2,-3){20}}
\put(140,130){\line(-1,-1){20}}  \put(260,130){\line(1,-1){20}}
\put(180,150){\line(-1,-1){60}}  \put(220,150){\line(1,-1){60}}
\put(240,130){\line(1,-2){30}}   \put(200,150){\line(-1,-2){30}}
\put(200,150){\line(1,-2){30}}   \put(200,110){\line(-1,-2){10}}
\put(200,110){\line(1,-2){10}}   \put(190,90){\line(1,-2){10}}
\put(210,90){\line(-1,-2){10}}   \put(170,90){\line(1,-2){30}}
\put(230,90){\line(-1,-2){30}}   \put(180,150){\line(1,-2){10}}
\put(220,150){\line(-1,-2){10}}  \put(30,100){\line(1,-1){40}}
\put(370,100){\line(-1,-1){40}}  \put(50,100){\line(1,-2){20}}
\put(350,100){\line(-1,-2){20}}  \put(90,100){\line(4,-3){40}}\qbezier(90,100)(93,100)(98,100)
\put(310,100){\line(-4,-3){40}}\qbezier(300,100)(303,100)(308,100)
\put(90,100){\line(-1,-2){10}}   \put(310,100){\line(1,-2){10}}
\put(80,80){\line(1,-1){10}}     \put(320,80){\line(-1,-1){10}}
\put(190,170){\line(-2,1){20}}   \put(210,170){\line(2,1){20}}
\multiput(160,200)(20,0){2}{\line(-1,-2){10}}
\multiput(220,200)(20,0){2}{\line(1,-2){10}}
\multiput(120,220)(10,-20){2}{\line(-3,-1){60}}
\multiput(280,220)(-10,-20){2}{\line(3,-1){60}}
\put(110,180){\line(-1,-1){20}}  \put(290,180){\line(1,-1){20}}
\put(90,210){\line(4,-1){40}}    \put(310,210){\line(-4,-1){40}}
\put(190,170){\line(1,-2){10}}   \put(210,170){\line(-1,-2){10}}
\put(130,70){\line(-6,-1){60}}   \put(270,70){\line(6,-1){60}}
\put(160,130){\line(-1,-2){30}}  
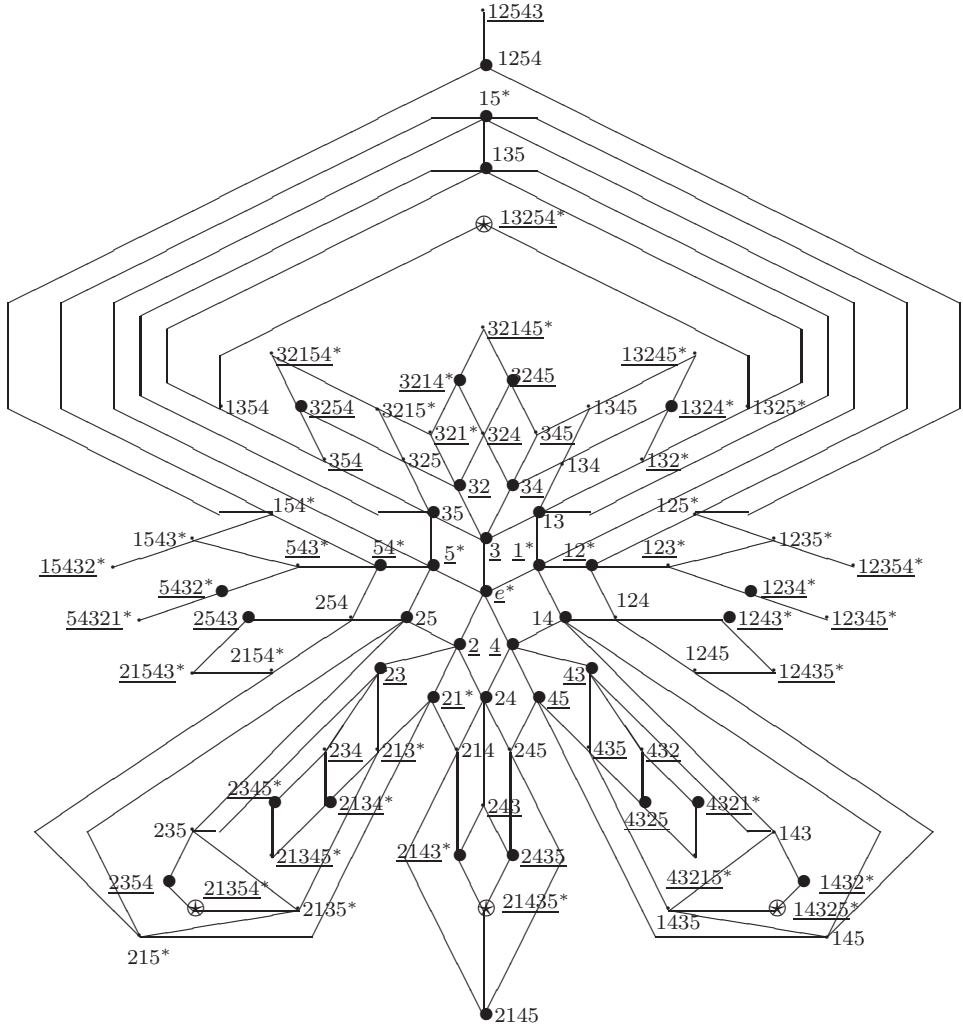
\put(140,0){\textbf{$\Gamma(bS_6)$} and \textbf{$\Gamma(bS_6)\bigcap\Gamma(*S_6)$}}
\end{picture}
\end{center}

\begin{remark}a) The figure contains the 89 braid-simple permutations in $\Gamma(\Sigma_6)$.

b) Among them, there are 58 elements in $bS_{6}\bigcap cS_{6}$ marked with $*$, 39 elements in $bS_{6}\bigcap gS_{6}$ marked with $\bullet$, 4 elements in $bS_{6}\bigcap sS_{6}$ marked with an inscribed $\star$, and 44 elements in $bS_{6}\bigcap tS_{6}$ which are underlined.
\end{remark}

\begin{center}
\begin{picture}(230,100)
\multiput(10,10)(100,0){3}{\line(0,1){80}}
\put(10,10){\line(5,4){100}}              \put(10,10){\line(5,2){200}}
\put(110,10){\line(5,4){100}}             \put(110,10){\line(-5,4){100}}
\put(210,10){\line(-5,4){100}}            \put(210,10){\line(-5,2){200}}
\put(8,8){\tiny$\bullet$\tiny$e$}         \put(8,88){\tiny$\bullet$\tiny$1$}
\put(108,8){\tiny$\bullet$\tiny$136$}     \put(108,88){\tiny$\bullet$\tiny$3$}
\put(208,8){\tiny$\bullet$\tiny$26$}      \put(208,88){\tiny$\bullet$\tiny$6$}
\put(38,64){\tiny$\bullet$\tiny$13$}      \put(110,23){\tiny$\bullet$\tiny$36$}
\put(188,71){\tiny$\bullet$\tiny$16$}     \put(60,68){\tiny$\bullet$\tiny$14$}
\put(90,56){\tiny$\bullet$\tiny$4$}       \put(130,40){\tiny$\bullet$\tiny$24$}
\put(180,20){\tiny$\bullet$}              \put(170,15){\tiny$246$}
\put(125,75){\tiny$\bullet$\tiny$35$}     \put(150,55){\tiny$\bullet$\tiny$5$}
\put(170,39){\tiny$\bullet$\tiny$25$}     \put(195,20){\tiny$\bullet$\tiny$2$}
\put(45,-15){\textbf{A $K_{3,3}$ subgraph of $\Gamma(bS_{7})$}}
\end{picture}
\end{center}

\bigskip

\bigskip

Direct computations give the next results:
\begin{lemma}\label{computation}Let $\pi$ be a cycle in $\Sigma_{n}$.

\noindent a) If $i+1$ is not an element of this cycle, then
$$\pi\tau_{i}=(k\ldots a, i,b\ldots \,h)(i,i+1)=(k\,\ldots\,a,i,i+1,b\ldots\,h).$$
\noindent b) If $i$ is not an element of this cycle, then
$$\pi\tau_{i}=(k\ldots a, i+1,b\ldots \,h)(i,i+1)=(k\,\ldots\,a,i+1,i,b\ldots\,h).$$
\noindent c) The next equalities hold
$$(k\ldots a,i, i\pm1,b\ldots \,h)(i,i\pm 1)=(k\,\ldots\,a,i,b\ldots\,h).$$
\noindent d) If $i,i+1$ are non consecutive elements of this cycle, then $\pi\tau_{i}$ is a product of two disjoint cycles:
\begin{eqnarray*}
  \pi\tau_{i} &=& (k\ldots a,i,b\ldots\,c,i+1,d\ldots \,h)(i,i+1) \\
   &=& (\ldots\,c,i+1,b\ldots)(k\,\ldots\,a,i,d\ldots\,h).
\end{eqnarray*}
\end{lemma}
\begin{cor}\label{cor:computation}If $\pi$ is a product of two disjoint cycles, one containing $i$ and the other containing $i+1$, then
\begin{eqnarray*}
 \pi\tau_{i}&=&(m\,\ldots\,c,i+1,b\ldots\,l)(k\,\ldots\,a,i,d\ldots\,h)(i,i+1) \\
 &=& (k\,\ldots a,i,b\ldots\,c,i+1,d,\ldots\,h).
\end{eqnarray*}
\end{cor}
Previous computations explain the next definitions, necessary to describe the connected components of $\Gamma(cS_{n})$. First we introduce, by two examples, an (oriented) 
\begin{figure}[hbt]
\begin{center}
\begin{pspicture}(0,3)(10,0)
\pscircle(1.5,1.5){1.5}
\rput(-0.25,1.5){$4$}\psdot(0,1.5)\psline[arrowsize=6pt]{->}(3.2,1.5)(1.5,1.5)
\rput(3.25,1.5){$1$}\psdot(3,1.5)       \rput(0.25,0.45){$5$}\psdot(0.35,0.55)
\rput(2.7,0.43){$6$}\psdot(2.6,0.55)    \rput(0.30,2.62){$3$}\psdot(0.45,2.55)
\rput(2.7,2.62){$2$}\psdot(2.47,2.6)    \pspolygon(3,1.5)(0,1.5)(2.47,2.6)(0.35,0.55)(2.6,0.55)
\pscircle(8.5,1.5){1.5}                 \rput(6.75,1.5){$4$}
\psdot(7,1.5)                           \psline[arrowsize=6pt]{->}(10,1.5)(8,0.75)
\rput(10.25,1.5){$1$}\psdot(10,1.5)     \rput(7.25,0.45){$5$}\psdot(7.35,0.55)
\rput(9.7,0.43){$6$}\psdot(9.5,0.37)    \rput(7.30,2.62){$3$}\psdot(7.45,2.55)
\rput(9.7,2.62){$2$}\psdot(9.47,2.6)    \pspolygon(10,1.5)(7.45,2.55)(9.5,0.37)(9.47,2.6)(7.35,0.55)
\rput(-1.5,1.5){$(61425)$}              \rput(5.5,1.5){$(63152)$}
\end{pspicture}
\end{center}
\end{figure}
\emph{polygonal representation} of a cycle: see the next two diagrams. A \emph{reduction move} of an oriented $k$-gon  ($k\geq 4$) consists in replacing a side $i\rightarrow j$ by the vertex $i$, and the side $j\rightarrow k$ by $i\rightarrow k$, if the following condition is fulfilled:
the interval $(\min\{i,j\},\max \{i,j\})$ does not contain another vertex of the polygon. For instance, the first pentagon, (61425) can be reduced in two steps to the triangle (514) or (614) or (615), and the second pentagon, (63152), is irreducible. Any polygon can be reduced to a unique irreducible type  ($a_{1},a_{2},\ldots,\,a_{s}$) or can be reduced to a triangle (this is not unique). An \emph{irreducible} type $(a_{1},a_{2},\ldots,a_{s})$ is a sequence of distinct integers (called vertices) in the interval $[1,n]$ such that $a_{1}=\max (a_{i})$ and, for any $i$, in the interval $[\min(a_{i}, a_{i+1}), \max (a_{i}, a_{i+1})]$ there is at least one other vertex $a_{j}$. We introduce the \emph{neighboring intervals}  $I^{+}(a_{i})$ and $I^{-}(a_{i})$ as follows: $I^{+}(a_{i})=[a_{i}+1,a_{j}-1]$, where $a_{j}$ is the smallest vertex greater than $a_{i}$ (in the special case of $a_{1}$, $I^{+}(a_{1})=[a_{1}+1,n]$), and $I^{-}(a_{i})=[a_{h}+1,a_{i}-1]$, where $a_{h}$ is greatest vertex smaller than $a_{i}$, (in the special case of $a_{m}=\min(a_{i})$, $I^{-}(a_{m})=[1,a_{m}-1]$). In the next picture there are only three non empty neighboring intervals: $I^{+}(3)=I^{-}(5)=\{4\}$ and $I^{+}(7)=\{8\}:$
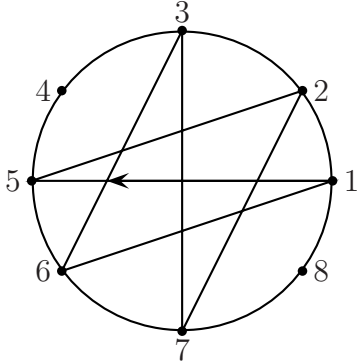
 \begin{figure}[hbt]
\begin{center}
\begin{pspicture}(0,4)(4,0)
 \pscircle(2,2){2}
 \psdot(0,2)\psdot(4,2)  \rput(-0.25,2){$5$}\rput(4.25,2){$1$}
 \psdot(2,0)\psdot(2,4) \rput(2,-0.25){$7$} \rput(2,4.25){$3$}
 \psdot(0.4,0.8)\psdot(3.6,0.8) \rput(0.15,0.8){$6$}\rput(3.85,0.8){$8$}
 \psdot(0.4,3.2)\psdot(3.6,3.2)\rput(0.15,3.2){$4$}\rput(3.85,3.2){$2$}
 \pspolygon(4,2)(0,2)(3.6,3.2)(2,0)(2,4)(0.4,0.8)
 \psline[arrowsize=6pt]{->}(4,2)(1,2)
\end{pspicture}
\end{center}
\end{figure}

Now it is easy to see that a polygon $P$ can be reduced to the irreducible type $(a_{1},\ldots,a_{s})$ if and only if it has the next structure:
$$ (a_{1},b^{1}_{1},\ldots,b^{1}_{t_{1}},a_{2},b_{1}^{2},\ldots,b^{2}_{t_{2}},\ldots,a_{s},b^{s}_{1},
\ldots,b^{s}_{t_{s}}), $$
where, for any $i$, all the elements $\{b_{1}^{i},\ldots,b_{t_{i}}^{i}\}$ are either in $I^{+}(a_{i})$ or in $I^{-}(a_{i})$ (in the case of equality $I^{+}(a_{i})=I^{-}(a_{j})$, the sets $\{a_{i},b^{i}_{*}\}$ and $\{a_{j},b_{*}^{j}\}$ should be separated: $\max b_*^i<\min b_*^j)$ and also the polygons $(a_{1},a_{2},\ldots,a_{i},b_{1}^{i},b_{2}^{i},\ldots,b_{t_{i}}^{i},a_{i+1},\ldots,a_{s})$ can be reduced to the irreducible type $(a_{1},\ldots,a_{i},a_{i+1},\ldots,a_{s})$. The uniqueness of irreducible types $(a_{1}\ldots a_{s})$ ($s\geq 5$) comes from the invariance of the unremovable vertices (the leftmost vertices $a$ in the previous formula).

\begin{proposition}\label{cycle graph is connected}a) The cycle-simple graph $\Gamma(cS_{n})$ is connected if and only if $n\leq 4$.

\noindent b) The connected component of identity in $\Gamma(cS_{n})$ contains only identity, all the transpositions, and all the cycles reducible to a triangle.

\noindent c) Any other component contains all the oriented polygons reducible to a given irreducible type $(a_{1},a_{2},\ldots,a_{s})$.
\end{proposition}
\begin{proof} a) For $n\leq 3$ every permutation is cycle-simple: $cS_{n}=\Sigma_{n}$ and for $n=4$ the graph $\Gamma(cS_{4})$ can be
seen in the diagram of $P(cS_4)$.. But if $n\geq 5$, using c),  we have at least two pentagonal types, (52413) and (53142):

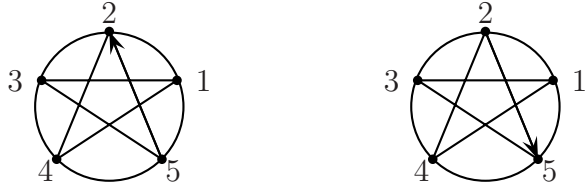
\begin{figure}[hbt]
\begin{center}
\begin{pspicture}(0,2.5)(8,0)
\pscircle(1,1){1}
\psdot(0.3,0.3)\psdot(1.7,0.3)\rput(0.15,0.15){$4$} \rput(1.85,0.15){$5$}
\psdot(1,2) \rput (1,2.25){$2$}
\psdot(0.1,1.35)\psdot(1.9,1.35)\rput(-0.25,1.35){$3$} \rput(2.25,1.35){$1$}
\pspolygon(0.3,0.3)(1.9,1.35)(0.1,1.35)(1.7,0.3)(1,2)
\psline[arrowsize=6pt]{->}(1.7,0.3)(1,2)
\rput(5,0){\pscircle(1,1){1}
\psdot(0.3,0.3)\psdot(1.7,0.3)\rput(0.15,0.15){$4$} \rput(1.85,0.15){$5$}
\psdot(1,2) \rput (1,2.25){$2$}
\psdot(0.1,1.35)\psdot(1.9,1.35)\rput(-0.25,1.35){$3$} \rput(2.25,1.35){$1$}
\pspolygon(0.3,0.3)(1.9,1.35)(0.1,1.35)(1.7,0.3)(1,2)
\psline[arrowsize=6pt]{->}(1,2)(1.7,0.3)}
\end{pspicture}
\end{center}
\end{figure}

b) We have to show that two cycles belong to the same connected component of $\Gamma(cS_{n})$  if and only if they have polygons related by reduction moves.  Lemma \ref{computation} shows that if two cycles $\gamma_{1}$, $\gamma_{2}$ are connected by $\tau_{i}$, $\gamma_{1}=\gamma_{2}\tau_{i}$  then one of them, say $\gamma_{2}$, has the form $\gamma_{2}=(k \ldots a,i,i+1,b \ldots h)$ or $\gamma_{2}=(k \ldots a,i+1,i,b \ldots h)$, hence the polygon $P(\gamma_{2})$  can be reduced to $P(\gamma_{1})$: remove the side $i\rightarrow i+1$ or $i+1\rightarrow i$, respectively.

Conversely, a polygonal reduction move 
$$P(\gamma_{2})=(\ldots i\rightarrow j \rightarrow k \ldots)\mapsto P(\gamma_{1})=(\ldots i\rightarrow k \ldots )$$ (we suppose that $i< j$) gives the following path in $\Gamma(cS_{n})$:

 $\gamma_{2}=(\ldots i,j,k \ldots)$---$(\ldots i,i+1,j,k\ldots)$---$\ldots$---$ (\ldots i,i+1\ldots j-1,j,k \ldots )$---

 ---$(\ldots i,i+1\ldots j-1,k\ldots)$---$\ldots$---$(\ldots i,i+1,k \ldots )$---$(\ldots i,k \ldots)=\gamma_{1}$.

\noindent The case $(\ldots i\rightarrow j\rightarrow k\ldots)$ with $i>j$ can be treated in the same way.
\end{proof}
\begin{example} $\Gamma(cS_{5})$ has three connected components: the connected component of identity and two isolated points: $(5\,2\,4\,1\,3)$ and $(5\,3\,1\,4\,2)$.
 \end{example}
 
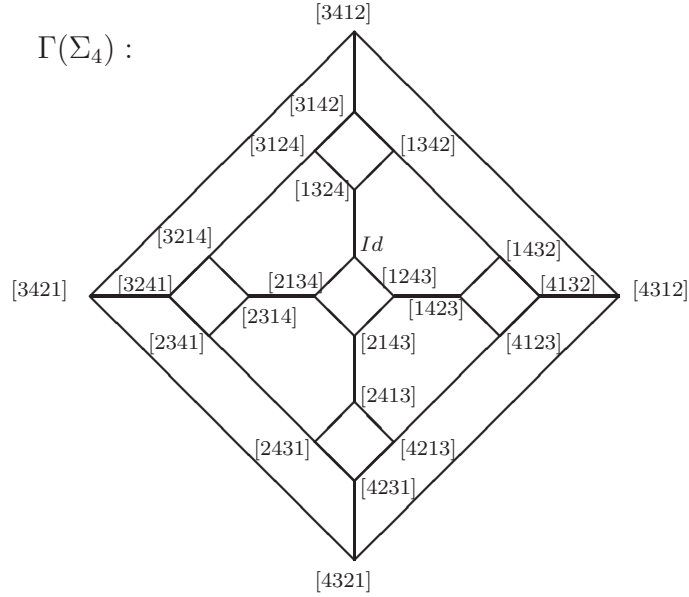
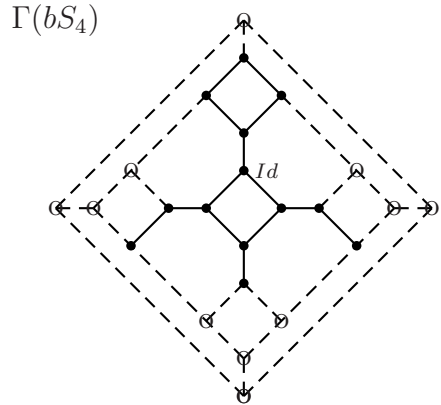
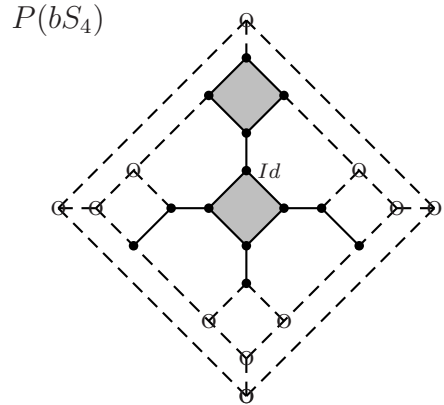
\begin{center}
\begin{picture}(220,220)
\thicklines
\put(100,0){\line(-1,1){100}} \put(100,30){\line(-1,1){70}} \put(100,0){\line(0,1){30}}
\put(100,0){\line(1,1){100}}  \put(100,30){\line(1,1){70}}  \put(0,100){\line(1,0){30}}
\put(200,100){\line(-1,1){100}} \put(170,100){\line(-1,1){70}}  \put(200,100){\line(-1,0){30}}
\put(0,100){\line(1,1){100}}      \put(30,100){\line(1,1){70}}  \put(100,200){\line(0,-1){30}}
\put(85,30){\SQR}    \put(85,140){\SQR}   \put(30,85){\SQR}
\put(85,85){\SQR}    \put(140,85){\SQR}   \put(85,100){\line(-1,0){25}}
\put(115,100){\line(1,0){25}}  \put(100,85){\line(0,-1){25}}  \put(100,115){\line(0,1){25}}
\put(85,-10){\tiny$[4321]$}   \put(205,100){\tiny$[4312]$}  \put(85,205){\tiny$[3412]$}
\put(-30,100){\tiny$[3421]$} \put(102,25){\tiny$[4231]$}    \put(117,40){\tiny$[4213]$}
\put(157,80){\tiny$[4123]$}   \put(170,102){\tiny$[4132]$}  \put(157,115){\tiny$[1432]$}
\put(117,155){\tiny$[1342]$}  \put(75,170){\tiny$[3142]$}   \put(60,155){\tiny$[3124]$}
\put(25,120){\tiny$[3214]$}  \put(10,102){\tiny$[3241]$}  \put(22,80){\tiny$[2341]$}
\put(62,40){\tiny$[2431]$}    \put(102,60){\tiny$[2413]$}  \put(102,80){\tiny$[2143]$}
\put(120,93){\tiny$[1423]$}  \put(110,105){\tiny$[1243]$} \put(102,117){\tiny$Id$}
\put(77,138){\tiny$[1324]$}   \put(67,103){\tiny$[2134]$}  \put(57,90){\tiny$[2314]$}
\put(-20,190){$\Gamma(\Sigma_{4}):$}
\end{picture}
\end{center}

\begin{figure}[hbt]
\begin{center}
\begin{pspicture}(0,6)(12,0)
\psline[linestyle=dashed](0,2.5)(2.5,0)   \psline[linestyle=dashed](2.5,0)(5,2.5)
\psline[linestyle=dashed](0,2.5)(2.5,5)   \psline[linestyle=dashed](5,2.5)(2.5,5)
\psline[linestyle=dashed](2.5,0.5)(0.5,2.5)  \psline[linestyle=dashed](2.5,0.5)(4.5,2.5)
\psline[linestyle=dashed](0.5,2.5)(2,4)    \psline[linestyle=dashed](4.5,2.5)(3,4)
\psline[linestyle=dashed](2,1)(2.5,1.5)    \psline[linestyle=dashed](3,1)(2.5,1.5)
\psline[linestyle=dashed](1.5,2.5)(1,3)    \psline[linestyle=dashed](3.5,2.5)(4,3)
\psline[linestyle=dashed](2.5,0)(2.5,0.5)  \psline[linestyle=dashed](2.5,4.5)(2.5,5)
\psline[linestyle=dashed](0,2.5)(0.5,2.5)  \psline[linestyle=dashed](4.5,2.5)(5,2.5)
\pspolygon(2,4)(2.5,3.5)(3,4)(2.5,4.5)      \pspolygon(2,2.5)(2.5,2)(3,2.5)(2.5,3)
\psline(2.5,1.5)(2.5,2)  \psline(2.5,3)(2.5,3.5)  \psline(1,2)(1.5,2.5)
\psline(1.5,2.5)(2,2.5)  \psline(3,2.5)(3.5,2.5)  \psline(3.5,2.5)(4,2)
\psdot*(1,2)  \psdot*(1.5,2.5)  \psdot*(2,2.5)
\psdot*(2.5,2)  \psdot*(2.5,1.5)  \psdot*(3,2.5)
\psdot*(3.5,2.5) \psdot*(4,2) \psdot*(2.5,3)
\psdot*(2.5,3.5)\psdot*(2.5,4.5)\psdot*(2,4)\psdot*(3,4)
\rput(0,2.5){$\text{o}$}\rput(2.8,3){\tiny$Id$}
\rput(2.5,0){$\text{o}$} \rput(5,2.5){$\text{o}$}
\rput(2.5,5){$\text{o}$}  \rput(0.5,2.5){$\text{o}$}
\rput(1,3){$\text{o}$}\rput(4,3){$\text{o}$}  \rput(2.5,0.5){$\text{o}$}
\rput(4.5,2.5){$\text{o}$}\rput(3,1){$\text{o}$}
\rput(2.5,0,5){$\text{o}$}\rput(2,1){$\text{o}$}
\rput(0,5){$\Gamma(bS_{4})$}
\psline[linestyle=dashed](7,2.5)(9.5,0)   \psline[linestyle=dashed](9.5,0)(12,2.5)
\psline[linestyle=dashed](7,2.5)(9.5,5)   \psline[linestyle=dashed](12,2.5)(9.5,5)
\psline[linestyle=dashed](9.5,0.5)(7.5,2.5)  \psline[linestyle=dashed](9.5,0.5)(11.5,2.5)
\psline[linestyle=dashed](7.5,2.5)(9,4)   \psline[linestyle=dashed](11.5,2.5)(10,4)
\psline[linestyle=dashed](9,1)(9.5,1.5)   \psline[linestyle=dashed](10,1)(9.5,1.5)
\psline[linestyle=dashed](8.5,2.5)(8,3)   \psline[linestyle=dashed](10.5,2.5)(11,3)
\psline[linestyle=dashed](9.5,0)(9.5,0.5)  \psline[linestyle=dashed](9.5,4.5)(9.5,5)
\psline[linestyle=dashed](7,2.5)(7.5,2.5)  \psline[linestyle=dashed](11.5,2.5)(12,2.5)
\pspolygon[linecolor=lightgray,fillstyle=solid, fillcolor=lightgray](9,4)(9.5,3.5)(10,4)(9.5,4.5)
\pspolygon[linecolor=lightgray,fillstyle=solid, fillcolor=lightgray](9,2.5)(9.5,2)(10,2.5)(9.5,3)
\pspolygon(9,4)(9.5,3.5)(10,4)(9.5,4.5)  \pspolygon(9,2.5)(9.5,2)(10,2.5)(9.5,3)
\psline(9.5,1.5)(9.5,2) \psline(9.5,3)(9.5,3.5) \psline(8,2)(8.5,2.5) 
\psline(8.5,2.5)(9,2.5) \psline(10,2.5)(10.5,2.5) \psline(10.5,2.5)(11,2)
\psdot*(8,2)\psdot*(8.5,2.5)\psdot*(9,2.5)
\psdot*(9.5,2)\psdot*(9.5,1.5)\psdot*(10,2.5)
\psdot*(10.5,2.5)\psdot*(11,2)\psdot*(9.5,3)  \rput(9.5,0.5){$\text{o}$}
\psdot*(9.5,3.5)\psdot*(9.5,4.5)\psdot*(9,4)\psdot*(10,4)
\rput(7,2.5){$\text{o}$}\rput(9.8,3){\tiny$Id$} \rput(9.5,0){$\text{o}$}
\rput(12,2.5){$\text{o}$}   \rput(9.5,5){$\text{o}$}
\rput(7.5,2.5){$\text{o}$}  \rput(8,3){$\text{o}$}\rput(11,3){$\text{o}$}
\rput(11.5,2.5){$\text{o}$}\rput(10,1){$\text{o}$}
\rput(9.5,0,5){$\text{o}$}\rput(9,1){$\text{o}$}
\rput(7,5){$P(bS_{4})$}
\end{pspicture}
\end{center}
\end{figure}

\begin{figure}[hbt]
\begin{center}
\begin{pspicture}(0,6)(12,0)
\psline[linestyle=dashed](0,2.5)(2.5,0)    \psline[linestyle=dashed](2.5,0)(5,2.5)
\psline[linestyle=dashed](0,2.5)(2.5,5)    \psline[linestyle=dashed](5,2.5)(2.5,5)
\psline[linestyle=dashed](2.5,0.5)(4.5,2.5)   \psline[linestyle=dashed](2.5,0.5)(0.5,2.5)
\psline[linestyle=dashed](2.5,4.5)(0.5,2.5)   \psline[linestyle=dashed](2.5,4.5)(4.5,2.5)
\psline[linestyle=dashed](2,1)(2.5,1.5)   \psline[linestyle=dashed](3,1)(2.5,1.5)
\psline[linestyle=dashed](1,2)(1.5,2.5)   \psline[linestyle=dashed](1,3)(1.5,2.5)
\psline[linestyle=dashed](2,4)(2.5,3.5)   \psline[linestyle=dashed](3,4)(2.5,3.5)
\psline[linestyle=dashed](4,2)(3.5,2.5)   \psline[linestyle=dashed](4,3)(3.5,2.5)
\psline[linestyle=dashed](2.5,1.5)(2.5,2)  \psline[linestyle=dashed](2.5,3)(2.5,3.5)
\psline[linestyle=dashed](1.5,2.5)(2,2.5)   \psline[linestyle=dashed](3.5,2.5)(3,2.5)
\psline[linestyle=dashed](2.5,2)(3,2.5)   \psline[linestyle=dashed](2.5,2)(2,2.5)
\psline[linestyle=dashed](2.5,3)(2,2.5)   \psline[linestyle=dashed](2.5,3)(3,2.5)
\psline[linestyle=dashed](2.5,0)(2.5,0.5)  \psline[linestyle=dashed](2.5,4.5)(2.5,5)
\psline[linestyle=dashed](0,2.5)(0.5,2.5)  \psline[linestyle=dashed](4.5,2.5)(5,2.5)
\rput(1,2){$\text{o}$}\rput(1.5,2.5){$\text{o}$}  \rput(2,2.5){$\text{o}$}\rput(2.5,2){$\text{o}$}
\rput(3,2.5){$\text{o}$}\rput(3.5,2.5){$\text{o}$}  \rput(4,2){$\text{o}$}\rput(2.5,3){$\text{o}$}
\rput(2.5,3.5){$\text{o}$}\rput(2,4){$\text{o}$}  \rput(3,4){$\text{o}$}
\rput(0,2.5){$\text{o}$}   \rput(2.5,0){$\text{o}$}  \rput(5,2.5){$\text{o}$}
\rput(2.5,5){$\text{o}$}   \rput(0.5,2.5){$\text{o}$}
\rput(1,3){$\text{o}$}\rput(4,3){$\text{o}$}
\rput(4.5,2.5){$\text{o}$}\rput(3,1){$\text{o}$}
\rput(2.5,0.5){$\text{o}$}\rput(2,1){$\text{o}$}\rput(2.8,3){\tiny$Id$}
\psdot*(2.5,4.5)  \psdot*(2.5,1.5)  \rput(0,5){$P(sS_{4})$}
\psline[linestyle=dashed](7,2.5)(9.5,0)  \psline[linestyle=dashed](9.5,0)(12,2.5)
\psline[linestyle=dashed](12,2.5)(9.5,5)  \psline[linestyle=dashed](11.5,2.5)(9.5,4.5)
\psline[linestyle=dashed](9,4)(9.5,3.5)   \psline[linestyle=dashed](10,4)(9.5,3.5)
\psline[linestyle=dashed](9.5,3.5)(9.5,3)  \psline[linestyle=dashed](9.5,0)(9.5,0.5)
\psline[linestyle=dashed](11.5,2.5)(12,2.5)
\pspolygon[linecolor=lightgray,fillstyle=solid, fillcolor=lightgray](7,2.5)(9.5,5)(9.5,4.5)(7.5,2.5)
\pspolygon[linecolor=lightgray,fillstyle=solid, fillcolor=lightgray](7.5,2.5)(8,2)(8.5,2.5)(8,3)
\pspolygon[linecolor=lightgray,fillstyle=solid, fillcolor=lightgray](8,2)(9,1)(9.5,1.5)(9.5,2)(9,2.5)(8.5,2.5)
\pspolygon[linecolor=lightgray,fillstyle=solid, fillcolor=lightgray](9,1)(9.5,0.5)(10,1)(9.5,1.5)
\pspolygon[linecolor=lightgray,fillstyle=solid, fillcolor=lightgray](9.5,1.5)(10,1)(11,2)(10.5,2.5)(10,2.5)(9.5,2)
\pspolygon[linecolor=lightgray,fillstyle=solid, fillcolor=lightgray](10.5,2.5)(11,2)(11.5,2.5)(11,3)
\pspolygon[linecolor=lightgray,fillstyle=solid, fillcolor=lightgray](9,2.5)(9.5,2)(10,2.5)(9.5,3)
\pspolygon(7,2.5)(9.5,5)(9.5,4.5)(7.5,2.5)  \pspolygon(7.5,2.5)(8,2)(8.5,2.5)(8,3)
\pspolygon(8,2)(9,1)(9.5,1.5)(9.5,2)(9,2.5)(8.5,2.5)  \pspolygon(9,1)(9.5,0.5)(10,1)(9.5,1.5)
\pspolygon(9.5,1.5)(10,1)(11,2)(10.5,2.5)(10,2.5)(9.5,2)  \pspolygon(10.5,2.5)(11,2)(11.5,2.5)(11,3)
\pspolygon(9,2.5)(9.5,2)(10,2.5)(9.5,3)  \psdot(8,2)\psdot(8.5,2.5)
\psdot(9,2.5)\psdot(9.5,2)  \psdot(10,2.5)\psdot(10.5,2.5)  \psdot(11,2)\psdot(9.5,3)
\psdot(9,4)   \psdot(7,2.5)    \psdot(9.5,5)  \psdot(7.5,2.5)
\psdot(8,3)\psdot(11,3)  \psdot(11.5,2.5)\psdot(10,1)  \psdot(9.5,0.5)\psdot(9,1)
\psdot*(9.5,4.5)  \psdot*(9.5,1.5)  \rput(9.5,0){$\text{o}$}  \rput(12,2.5){$\text{o}$}
\rput(9.5,3.5){$\text{o}$}\rput(9.8,3){\tiny$Id$}  \rput(10,4){$\text{o}$}
\rput(7,5){$P(tS_{4})$}
\end{pspicture}
\end{center}
\end{figure}

\begin{figure}[hbt]
\begin{center}
\begin{pspicture}(0,5)(12,0)
\psline[linestyle=dashed](0,2.5)(2.5,0)  \psline[linestyle=dashed](2.5,0)(5,2.5)
\psline[linestyle=dashed](0,2.5)(2.5,5)  \psline[linestyle=dashed](5,2.5)(2.5,5)
\pspolygon[linecolor=lightgray,fillstyle=solid, fillcolor=lightgray](2,1)(2.5,0.5)(3,1)(2.5,1.5)
\pspolygon[linecolor=lightgray,fillstyle=solid, fillcolor=lightgray](0.5,2.5)(1,2)(1.5,2.5)(1,3)
\pspolygon[linecolor=lightgray,fillstyle=solid, fillcolor=lightgray](3.5,2.5)(4,2)(4.5,2.5)(4,3)
\pspolygon[linecolor=lightgray,fillstyle=solid, fillcolor=lightgray](2,4)(2.5,3.5)(3,4)(2.5,4.5)
\pspolygon[linecolor=lightgray,fillstyle=solid, fillcolor=lightgray](3,2.5)(3.5,2.5)(4,3)(3,4)(2.5,3.5)(2.5,3)
\pspolygon[linecolor=lightgray,fillstyle=solid, fillcolor=lightgray](1,3)(1.5,2.5)(2,2.5)(2.5,3)(2.5,3.5)(2,4)
\pspolygon(2,1)(2.5,0.5)(3,1)(2.5,1.5)  \pspolygon(0.5,2.5)(1,2)(1.5,2.5)(1,3)
\pspolygon(3.5,2.5)(4,2)(4.5,2.5)(4,3)   \pspolygon(2,4)(2.5,3.5)(3,4)(2.5,4.5)
\pspolygon(3,2.5)(3.5,2.5)(4,3)(3,4)(2.5,3.5)(2.5,3) \pspolygon(1,3)(1.5,2.5)(2,2.5)(2.5,3)(2.5,3.5)(2,4)
\psline(2,1)(1,2)  \psline(3,1)(4,2)  \psline[linestyle=dashed](2.5,1.5)(2.5,2)
\psline[linestyle=dashed](2.5,2)(2,2.5)  \psline[linestyle=dashed](2.5,2)(3,2.5)
\psline[linestyle=dashed](2.5,0)(2.5,0.5)  \psline[linestyle=dashed](2.5,4.5)(2.5,5)
\psline(0,2.5)(0.5,2.5)  \psline(4.5,2.5)(5,2.5)
\psdot*(1,2) \psdot*(1.5,2.5) \psdot*(2,2.5) \psdot*(2.5,1.5) \psdot*(3,2.5) \psdot*(3.5,2.5)
\psdot*(4,2) \psdot*(2.5,3) \psdot*(2.5,3.5) \psdot*(2.5,4.5) \psdot*(2,4) \psdot*(3,4)
\psdot(0,2.5) \psdot(0.5,2.5) \psdot(5,2.5)  \psdot*(2.5,0.5)
\rput(2.5,5){$\text{o}$} \rput(2.8,3){\tiny$Id$}  \rput(2.5,0){$\text{o}$}
\rput(2.5,2){$\text{o}$}  \psdot(1,3)\psdot(4,3)  \psdot(4.5,2.5)\psdot(3,1)
\psdot(2,1)  \rput(0,5){$P(cS_{4})$}
\psline[linestyle=dashed](7,2.5)(9.5,0)  \psline[linestyle=dashed](9.5,0)(12,2.5)
\psline[linestyle=dashed](7,2.5)(9.5,5)  \psline[linestyle=dashed](12,2.5)(9.5,5)
\psline[linestyle=dashed](9,1)(7.5,2.5)  \psline[linestyle=dashed](10,1)(11.5,2.5)
\psline[linestyle=dashed](9,1)(9.5,1.5)  \psline[linestyle=dashed](9.5,1.5)(9.5,2)
\psline[linestyle=dashed](9.5,1.5)(10,1)  \psline[linestyle=dashed](9,4)(9.5,4.5)
\psline[linestyle=dashed](10,4)(9.5,4.5)  \psline[linestyle=dashed](9.5,4.5)(9.5,5)
\psline[linestyle=dashed](8,2)(8.5,2.5)   \psline[linestyle=dashed](10.5,2.5)(11,2)
\psline[linestyle=dashed](7,2.5)(7.5,2.5)  \psline[linestyle=dashed](11.5,2.5)(12,2.5)
\pspolygon[linecolor=lightgray,fillstyle=solid, fillcolor=lightgray](10,2.5)(10.5,2.5)(11,3)(10,4)(9.5,3.5)(9.5,3)
\pspolygon[linecolor=lightgray,fillstyle=solid, fillcolor=lightgray](8,3)(8.5,2.5)(9,2.5)(9.5,3)(9.5,3.5)(9,4)
\pspolygon[linecolor=lightgray,fillstyle=solid, fillcolor=lightgray](9,2.5)(9.5,2)(10,2.5)(9.5,3)
\pspolygon(10,2.5)(10.5,2.5)(11,3)(10,4)(9.5,3.5)(9.5,3)  \pspolygon(8,3)(8.5,2.5)(9,2.5)(9.5,3)(9.5,3.5)(9,4)
\pspolygon(9,2.5)(9.5,2)(10,2.5)(9.5,3)  \psline(7.5,2.5)(8,3)  \psline(11,3)(11.5,2.5)
\psline(9.5,0)(9.5,0.5)  \psline(9.5,0.5)(9,1)  \psline(9.5,0.5)(10,1)
\psdot*(8.5,2.5)\psdot*(9,2.5)\psdot*(10,2.5)\psdot*(10.5,2.5)  \psdot*(9.5,0.5)
\psdot*(9.5,3)\psdot*(9.5,3.5)\psdot*(9,4)\psdot*(10,4)\psdot(7.5,2.5)
\psdot(9.5,5)\psdot(9.5,0)\psdot(8,3)\psdot(11,3)  \psdot(11.5,2.5)\psdot(10,1)
\psdot(9.5,2)  \psdot(9,1)  \rput(8,2){$\text{o}$}  \rput(11,2){$\text{o}$}
\rput(9.5,1.5){$\text{o}$}\rput(9.8,3){\tiny$Id$}  \rput(7,2.5){$\text{o}$}
\rput(12,2.5){$\text{o}$}  \rput(9.5,4.5){$\text{o}$}  \rput(7,5){$P(gS_{4})$}
\end{pspicture}
\end{center}
\end{figure}

\begin{proposition}\label{group graph is connected}a) The graph $\Gamma(gS_{n})$ is connected if and only if $n\leq 3$.

\noindent b) The connected component of identity contains only identity, all the transpositions $\tau_{i}$ and $(i+2,i)$, products of disjoint transpositions of two types: $\tau_{k_{1}}\tau_{k_{2}}\ldots \tau_{k_{s}}$ and $(i+2,i)\tau_{k_{1}}\tau_{k_{2}}\ldots \tau_{k_{s}}$, and also three three-cycles of the form:
 $$(i+2,i+1,i),(i+2,i,i+1),(i+3,i,i+2),(i+3,i+1,i);$$
\noindent c) the other components are divided in three classes: $C(j,i), C(J_{*}, I_{*})$ and isolated components:

c1) a component of type $C(j,i)$, where $j\geq i+3$, contains only $(j\,i), (j\,i)\tau_{k_{1}}\tau_{k_{2}}\ldots \tau_{k_{s}}$ (disjoint transpositions) and also (at most four) three cycles $(j,i,i+1),(j,i,i-1),(j+1,i,j),(j,j-1,i)$;

c2) a component of type $C(J_{*},I_{*})$, where $J_{*}=(j_{1},\ldots,j_{A})$, $I_{*}=(i_{1},\ldots,i_{A})$, ($A\geq 2$), $j_{a}\geq i_{a}+2$ for any $a$, and the cycles $(j_{1},i_{1}),\ldots,(j_{A},i_{A})$ are disjoint,  contains only the products of disjoint transpositions $(j_{1},i_{1})\ldots(j_{A},i_{A})\tau_{k_{1}}\ldots \tau_{k_{s}}$;

c3) singletons  of type $C^{+}_{k\,j\,i}=(k,j,i)$ or $C^{+}_{k\,j\,i}=(k,i,j)$, where $k-2\geq j\geq i+2$  and products of  disjoint such cycles
$C^{\pm}_{k_{1}\,j_{1}\,i_{1}}\ldots C^{\pm}_{k_{s}\,j_{s}\,i_{s}}$, and also  singletons of type $C_{p,q}$, where $p\geq 5$ is prime and $pq\leq n$, containing a unique product of disjoint $q$ cycles of length $p$.
 \end{proposition}
\begin{proof} All these are consequences of Lemma \ref{computation} and of the computations of the length of the cycles in those formulae: in Corollary~\ref{cor:computation},  $\alpha\beta\tau_{i}=(\ldots i \ldots )(\ldots i+1 \ldots)(i+1,i)=\gamma$, we have 
length$(\gamma)=\,\,$length$(\alpha)+$length$(\beta)$, without solution in primes (all the lengths should be equal to a unique prime).

In the case a) of Lemma \ref{computation}, $\alpha\tau_{i}=(\ldots i \ldots )(i+1,i)=\widetilde{\alpha}$, we have length $(\widetilde{\alpha})=$length$(\alpha)\pm 1$, with   solutions the primes 2 and 3. These give isolated g-simple permutations for primes $\geq 5$ and connected components corresponding to the primes 2 and 3, and also isolated components containing products of cycles of length 3 of types $C^{\pm}_{k\,i\,j}$ with $k-2 \geq i \geq j+2$.
 \end{proof}

\noindent \emph{Proof of Theorem \ref{bSn is contractible}.}   The contractibility of  $B(bS_n)$ and of 
$W(bS_n)$ is is a consequence of a general theorem of Quillen \cite{Quillen}: the posets
 ($bS_{n}$, Bruhat order) and ($bS_{n}$, weak order) have Id as a smallest element.

The contractibility of $P(bS_{n})$, the braid-simple part of the permutahedron, is given by a double induction, 
on $n$ and on $j$, the index of the next decreasing filtration:
$$P(bS_{n})=F_{n+1}\subset F_{n}\subset F_{n-1}\subset \ldots \subset F_{2}\subset F_{1}=P(bS_{n+1}),$$
where the $j$ stage of the filtration is
$$F_{j}=\left\{
          \begin{array}{c}
            \text{faces in $P(bS_{n+1})$ with vertices from $F_{j+1}$ and} \\
            \text{new vertices $D(k_{1},j_{1})\ldots D(k_{s-1},j_{s-1})D(n,j)$} \\
          \end{array}
        \right\}.
$$
The induction on $n$ starts with $n=1,2$ where $P(bS_{n})=P(\Sigma_{n})$  are homeomorphic to a point (Id $\bullet$)
 and to a segment (Id $\bullet$---$\bullet\tau_{1}$). We show that $F_{j+1}$ is deformation retract of $F_{j}$. The space $F_{j}$ is obtained from $F_{j+1}$ by adding vertices of the form $D(k_{1},j_{1})\ldots D(n,j)$ along the one dimensional cells $D(k_{1},j_{1})\ldots D(n,j+1)$---$D(k_{1},j_{1})\ldots D(n,j)$. The cells in $F_{j}$ are the cells of $F_{j+1}$ and new cells of two types: \emph{cylinders} $C_{\alpha}\times I_{j}$, where $C_{\alpha}$ is cell in $F_{j+1}$ and $I_{j}$ corresponds to the edge $Id$--$\tau_{j}$, and the \emph{faces} of these cylinders. There is no cell of $F_{j}$  "parallel" to $F_{j+1}$, $C_{\alpha}\times I_{j}$, 
 
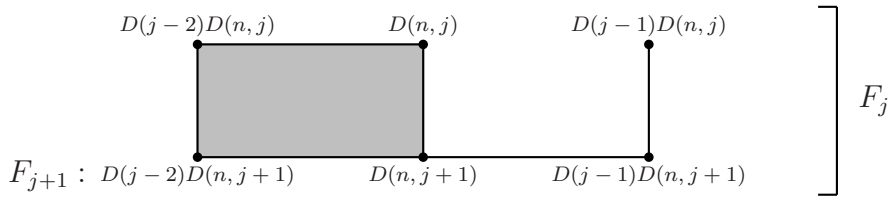
 \begin{figure}[hbt]
\begin{center}
\begin{pspicture}(0,3)(10,0)
\pspolygon[linecolor=lightgray,fillstyle=solid, fillcolor=lightgray](0.5,0.5)(3.5,0.5)(3.5,2)(0.5,2)
\pspolygon(0.5,0.5)(3.5,0.5)(3.5,2)(0.5,2)
\psline(3.5,0.5)(6.5,0.5)\psline(6.5,0.5)(6.5,2)
\psdot(0.5,0.5)\psdot(3.5,0.5)\psdot(3.5,2)\psdot(0.5,2)\psdot(6.5,0.5)\psdot(6.5,2)
\rput(0.5,0.25){\tiny$D(j-2)D(n,j+1)$}\rput(3.5,0.25){\tiny$D(n,j+1)$}\rput(6.5,0.25){\tiny$D(j-1)D(n,j+1)$}
\rput(0.5,2.25){\tiny$D(j-2)D(n,j)$}\rput(3.5,2.25){\tiny$D(n,j)$}\rput(6.5,2.25){\tiny$D(j-1)D(n,j)$}
\rput(-1.5,0.25){$F_{j+1}:$}
\psline(9,0)(9,2.5)\psline(9,0)(8.75,0)\psline(9,2.5)(8.75,2.5)
\rput(9.5,1.25){$F_{j}$}
\end{pspicture}
\end{center}
\end{figure}

\noindent if  $C_{\alpha}$ is not a cell in $F_{j}$: to a cell with vertices $\{D(K_{*}^{\lambda},J_{*}^{\lambda})D(n,j)\}_{\lambda \in \Lambda_{\alpha}}$ where $k_{1}^{\lambda}< k_{2}^{\lambda}< \ldots< k_{s-1}^{\lambda}\leq n-1$, corresponds a face $F_{\alpha}$ of $P(\Sigma_{n})$, hence $C_{\alpha}=F_{\alpha}\times (\tau_{n}\tau_{n-1}\ldots \tau_{j+1})$ is a face of $F_{j+1}$ (on the other hand, there are cells in $F_{j+1}$ with no parallel correspondent in $F_{j}$:)

 The deformation retract $\mathcal{D}_{t}:\,F_{j}\rightarrow F_{j}$, $\mathcal{D}_{1}=Id_{F_{j}}$, $\mathcal{D}_{0}:\,F_{j}\rightarrow F_{j+1}$ is defined by projecting the cylinders $C_{\alpha}\times I_{j}$ onto $C_{\alpha}$.    $\Box$

\end{document}